\documentclass[11pt]{article}
\usepackage{amscd,amssymb,stmaryrd}
\usepackage[fleqn]{amsmath}
\usepackage{graphicx}
\usepackage[utf8]{inputenc}
\usepackage[export]{adjustbox}
\usepackage{wrapfig}
\usepackage{amstext}
\usepackage{pstricks,pst-node,pst-plot,pst-coil}
\usepackage{amsthm}
\usepackage{bm}
\usepackage{amsmath}
\usepackage{color}
\usepackage{mathtools}
\usepackage{hyperref}
\usepackage[english]{babel}
\usepackage{booktabs}
\usepackage{mathrsfs}
\usepackage{multirow}
\usepackage{subfiles}
\usepackage{subfigure}
\usepackage{lineno}
\usepackage{floatrow}
\usepackage[T1]{fontenc}
\usepackage{caption}

\newtheorem{theorem}{Theorem}[section]

\newtheorem{lemma}[theorem]{Lemma}

\theoremstyle{definition}

\theoremstyle{remark}

\newfloatcommand{capbtabbox}{table}[][\FBwidth]

\title{Constraint energy minimization generalized multiscale finite element method in mixed formulation for parabolic equations}
\author{Yiran Wang\thanks{Department of Mathematics, The Chinese University of Hong Kong, Hong Kong SAR.}, \;
Eric Chung\thanks{Department of Mathematics, The Chinese University of Hong Kong, Hong Kong SAR.} \; and \;
Lina Zhao\thanks{Department of Mathematics, The Chinese University of Hong Kong, Hong Kong SAR.}
}
\begin{document}
\maketitle
\begin{abstract}
    In this paper, we develop the constraint energy minimization generalized multiscale finite element method (CEM-GMsFEM) in mixed formulation applied to parabolic equations with heterogeneous diffusion coefficients. The construction of the method is based on two multiscale spaces: pressure multiscale space and velocity multiscale space. The pressure space is constructed via a set of well-designed local spectral problems, which can be solved independently. Based on the computed pressure multiscale space, we will construct the velocity multiscale space by applying constrained energy minimization. The convergence of the proposed method is proved.
    In particular, we prove that the convergence of the method depends only on the coarse grid size, and is independent of the heterogeneities and contrast of the
    diffusion coefficient.
    Four typical types of permeability fields are exploited in the numerical simulations, and the results indicate that our proposed method works well
    and gives efficient and accurate numerical solutions.
\end{abstract}
\section{Introduction}
A great number of problems of fundamental and practical significance are described by partial differential equations with coefficients varying from a wide range of length scales. For example, flows in composite materials and porous media, transport in high Reynolds number flows are examples of this type. The heterogeneity and high-contrast properties of the coefficients cause extreme difficulty for the numerical approximation of these problems. In this paper, we consider a flow model described by a parabolic equation. As mentioned above, the multiscale nature of the problem should be taken into account. It is of great importance to develop effective approaches to solve these types of problems. Some examples of these approaches include upscaling methods \cite{chen2003coupled,durlofsky1991numerical,weh02,chung2018non} and multiscale methods \cite{efendiev2011multiscale,wheeler2012multiscale,efendiev2009multiscale,hou1997multiscale,jenny2003multi}.

In upscaling method, one typically upscales the media properties and solves the corresponding global problem on a coarse grid. This method may encounter inaccuracy when fine-grid information has strong effect on the solution. For multiscale methods, a significant convenience of this class of methods are the independent construction of a set of multiscale basis functions that span an approximation space associated with a coarse grid, whose grid size is much larger than the characteristic scale of the heterogeneous coefficient. A key idea in multiscale methods is that we construct multiscale basis functions using a set of local problems taking account of fine-scale information of media.
%
%
In practical applications, the mass conservation property for velocity fields are important, in line with this mechanism a variety of approaches are pursued such as multiscale finite volume methods \cite{cortinovis2014iterative,lunati2004multi}, mixed multiscale finite element methods \cite{aarnes2004use,aarnes2008mixed,ch03,chung2015mixed}, mortar multiscale methods \cite{arbogast2007multiscale, peszynska2005mortar,peszynska2002mortar}, post postprocessing methods \cite{odsaeter2017postprocessing,bush2013application}, etc. For example, in mixed multiscale finite element method, we utilize the mixed finite element formulation, where one may use first-order system for pressure and velocity. The pressure space is composed of piecewise constant functions. The support for each pressure basis function is a single coarse block while for velocity basis function, it vanishes outside a coarse neighborhood, which is composed of two coarse blocks sharing a coarse edge.

Furthermore, in cases where long channels and non-separable scales exist, only one basis function in each local neighborhood is not sufficient to guarantee a promising approximation. Therefore, one may consider taking multiple basis functions instead, which motivates the Generalized Multiscale Finite element method (GMsFEM)\cite{chung2016adaptiveJCP,efendiev2013generalized,chung2015perforated,WaveGMsFEM,eglp13,chan2016adaptive,chung2016mixed,chung2016mixedwave,chung2019mixedelastic,chung2017goal}. GMsFEM is a flexible general framework that generalizes the Multiscale Finite Element Method (MsFEM) by systematically enriching the coarse spaces. Before constructing offline space, where we take no consideration of some global information like source terms and boundary conditions, one needs to construct a snapshot space and use some well-designed local spectral problems to obtain the offline multiscale space. The construction of the local problem is motivated by the convergence analysis, which can give a convergence rate $1/\Lambda$, where $\Lambda$ is the smallest eigenvalue whose modes are excluded in the multiscale space. However, the convergence rate is independent of the coarse-mesh size $H$ while we show more interest in coarse-mesh-dependent convergence, which motivates the method we use in this paper.

In this paper, we consider a novel mass-conservative method to construct multiscale space, i.e. the constraint energy minimization generalized multiscale finite element method \cite{chung2018constraint,fu2020constraint} in mixed formulation (mixed CEM-GMsFEM) \cite{chung2018constraint_mixed,chung2020computational,chung2019online}, where the convergence rate of this method is proportional to $H/\Lambda$ under the condition that the number of oversampling layers is chosen suitably. The construction of the multiscale space is divided into two steps. First we construct an auxiliary space for pressure multiscale space by solving a set of well-designed eigenvalue problems in each coarse element. Then we compute velocity multiscale basis functions in the corresponding oversampled regions via energy minimization. We emphasize that by using the oversampling strategy and energy minimization,
most dominant features are localized within the oversampling region.
We prove the convergence rate $H/\Lambda$ under certain conditions on the oversampled layers. In particular, the size of oversampling domain depends weakly on the contrast of the media and based on this condition, the convergence rate is independent of the contrast.

In the numerical simulation, we consider four distinctive permeability fields to estimate the effects of the proposed method. In particular, we compare the effects of dimension of offline multiscale space, the number of oversampling layers and the coarse-mesh size separately. As for the permeability fields, minimal value is 1 in the first three cases while 1 is the largest value in the last case. As to the third case, the field is extracted from SPE 10 Benchmark tests \cite{christie2001tenth} while for the rest, the permeability fields have large contrast and distinct channels. In particular, the number of channels are selected to have a certain number of minimum basis functions in each coarse block. In terms of the third case, the number of channels are not explicitly defined. Our numerical results show a first-order convergence with respect to the coarse-mesh size. Consequently, one can achieve relatively low error by using well-chosen parameters, i.e. dimension of offline multiscale space, number of oversampling layers and the coarse-mesh size.

The paper is organized in the following way. In Section 2, we exhibit the concerned problem setting and provide some preliminaries. We will review the construction of multiscale space in Section 3. In Section 4, we will give the proof of the stability of the proposed method and convergence analysis. We will show the numerical results in Section 5. Finally, a concluding remark is given.
\section{Model formulation}

In this section we describe the model problem and some basic notations that will be exploited throughout the paper. Let $\Omega\subset \mathbb{R}^2$ be a bounded computational domain and we denote
$$ V_0:=\{v\in H(\text{div},\Omega):v\cdot \mathbf{n}=0 \quad\text{on }\partial\Omega\},\quad Q:=L^2(\Omega). $$
We seek the velocity $\mathbf{v}\in V_0$ and the pressure $p\in Q$ such that
\begin{eqnarray}
\begin{aligned}
\kappa^{-1} \mathbf{v}+\nabla{p} &=0 & & \text { in } \Omega \times(0, T], \\
\rho \dot{p}+\nabla \cdot \mathbf{v} &=f & & \text { in } \Omega \times(0, T], \\
\mathbf{v} \cdot \mathbf{n} &=0 & & \text { on } \partial \Omega \times[0, T], \\
\mathbf{v}(0,x) &=h_{\mathbf{v}} & & \text { in } \Omega, \\
p(0,x) &=h_{p} & & \text { in } \Omega.
\label{model}
\end{aligned}
\end{eqnarray}
Here we use $\dot{p}$ to denote the time derivative of $p$ and $T>0$ is a given terminal time.
The function $\rho\in L^{\infty}(\Omega)$ is the (positive) density of the fluid satisfying $0<\rho_{\text{min}}\leq \rho$ and as usual $\mathbf{n}$ is the unit outward normal vector to the boundary $\partial \Omega$. Here, the permeability field $\kappa:\Omega\rightarrow \mathbb{R}$ is assumed to be highly oscillatory and of high contrast. In particular, almost every $x\in \Omega$ satisfies  $\frac{\mathcal{K}_{\text{max}}}{\mathcal{K}_{\text{min}}}\gg 1$ combined with $0<\kappa_{\text{min}}\leq \kappa(x)\leq \kappa_{\text{max}}$. The source function satisfies $f\in L^2(\Omega)$. Here, $h_{\mathbf{v}}$ and $h_p$ are some given initial conditions.

Multiplying (\ref{model}) with corresponding test functions and performing integration by parts, we can obtain the following mixed weak formulation: Finding $\mathbf{v}\in V_0$ and $p\in Q$ such that
\begin{eqnarray}
\begin{aligned}
a(\mathbf{v},\mathbf{w})-b(\mathbf{w},p) &= 0\quad\forall \mathbf{w}\in V_0,\\
(\dot{p},q)_{\rho}+b(\mathbf{v},q)&=(f,q)\quad\forall q\in Q,\label{variational}
\end{aligned}
\end{eqnarray}
where the bilinear forms $a:V_0\times V_0\rightarrow \mathbb{R}$ and $b:V_0\times Q\rightarrow\mathbb{R}$ are defined as follows:
\begin{eqnarray*}
a(\mathbf{v},\mathbf{w}):=\int_{\Omega}\kappa^{-1}\mathbf{v}\cdot \mathbf{w}\; dx, \quad b(\mathbf{v},p):=\int_{\Omega}p\nabla\cdot \mathbf{v}\; dx
\end{eqnarray*}
for all $\mathbf{v},\mathbf{w}\in V_0$, and $p\in Q$. Here for $D\subset \mathbb{R}^2$, we denote the inner product in $L^2(\Omega)$ by $(\xi,\zeta)_D:=\int_D \xi \zeta\;dx$ for $\xi,\zeta\in L^2(\Omega)$, and we use $(\xi,\zeta)_{\rho}:=\int_D \rho \xi \zeta\;dx$ to denote the weighted inner product. When $D$ coincides with $\Omega$, the subscript is omitted.

We remark that $b(\cdot,\cdot)$ should satisfy the following inf-sup condition \cite{bathe2001inf}:
For all $q\in Q$, with $\int_{\Omega} q dx=0$, there exists a constant $C_0>0$ independent of $\kappa$ such that
\begin{eqnarray}
\|q\|_{L^2(\Omega)}\leq C_0 \sup_{\mathbf{v}\in V_0}\dfrac{b(\mathbf{v},q)}{\|\mathbf{v}\|_{H(div,\Omega)}}.\label{inf-sup}
\end{eqnarray}

To speed up the computation, we will then construct the constraint energy minimizing generalized multiscale finite element method (CEM-GMsFEM) in mixed formulation, which originated from GMsFEM. The key idea of this method is constructing local basis functions with minimal energy norm based on a set of well-designed local spectral problems. The details will be presented in the next section.

\begin{figure}[ht]
    \centering
    \subfigure
       { \includegraphics[width=2.5in]{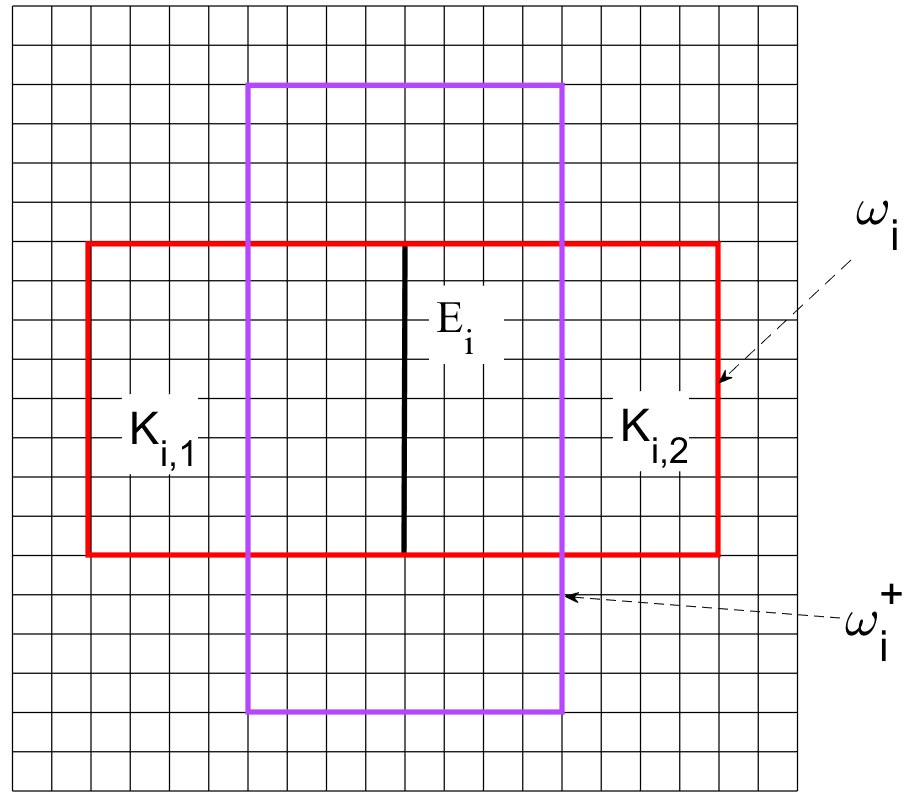}}
        \subfigure
       { \includegraphics[width=3.5in]{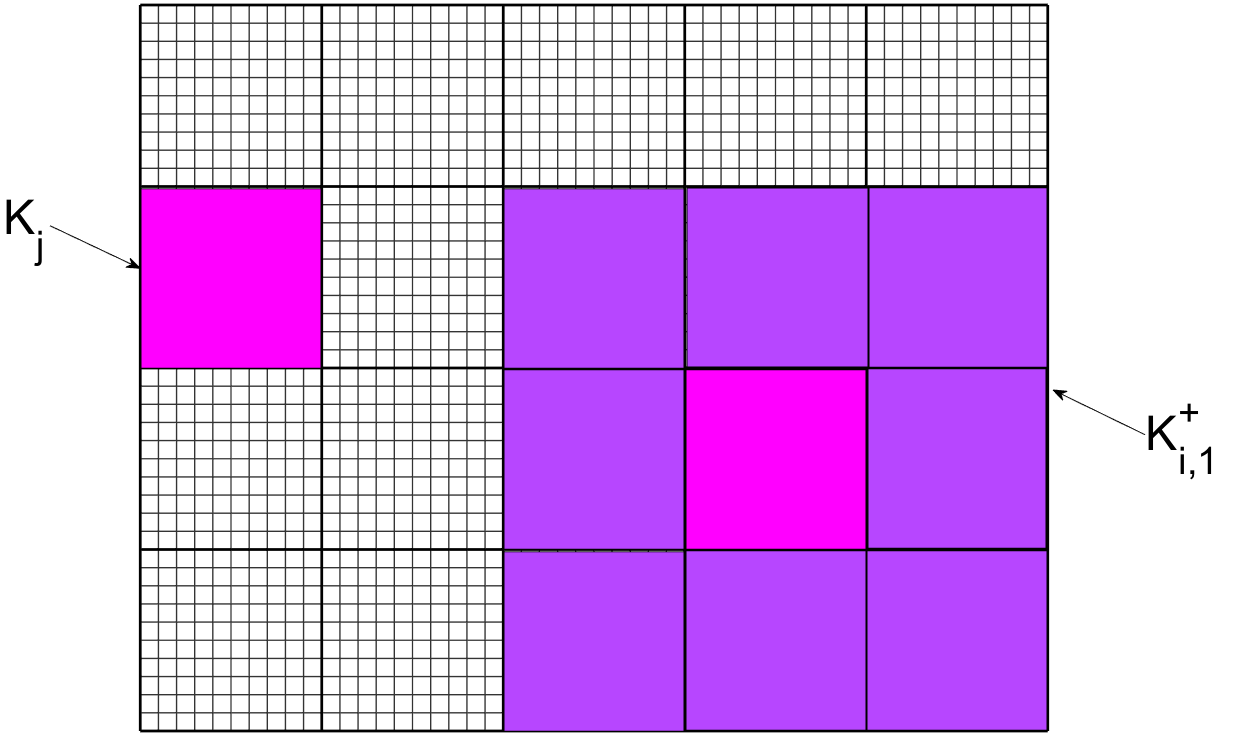}}
    \caption{Left: Coarse neighborhood and oversampled coarse neighborhood. Right: Coarse element and oversampled coarse element}\label{fig:mesh}
\end{figure}

Before closing this section we introduce some notations that will be used later. We consider a conforming rectangular partition $\mathcal{T}_H$ of the domain $\Omega$ with mesh size $H$ (defined in (\ref{H})). Under $\mathcal{T}_H$, $\Omega$ is a union of non-overlapping (excluding edges) coarse elements $K_i$, for $i=1,\ldots,N_e$, where $N_e$ is the total number of the coarse elements.
 We have
 \begin{eqnarray}
   H=\max_{i=1,\ldots,N_e}\max_{x,y\in K_i}|x-y|.\label{H}
 \end{eqnarray}
 Based on the coarse partition, we further define a finer partition $\mathcal{T}_h$ such that each coarse element is a union of disjoint fine elements. More specifically, $\forall K_{j} \in \mathcal{T}_{H}$, $ K_{j}=\bigcup_{F\in I_{j} }F$ for some $I_{j}\in \mathcal{T}_{h}$. In this paper, we consider identical rectangular partition and the method can be extended to general partition.  We denote the interior coarse edges of $\mathcal{T}_{H}$ by $E_i,i=1,\cdots,N_{\text{in}}$,
 where $N_\text{in}$ is the number of interior coarse edges. We define the coarse neighborhood of the edges $E_i$ by $\omega_i:=\cup\{K_j\in \mathcal{T}_{H}:E_i\subset \overline{K_j}\}$.
 An illustration of the mesh notations is shown in Figure \ref{fig:mesh}.
 We remark that each coarse neighborhood is composed of two coarse blocks sharing a common edge. As in Figure \ref{fig:mesh}, $\omega_i=K_{i,1}\cup K_{i,2}$. Besides, we define the oversampled coarse element $K_{i,l}^{+}$ by enlarging $K_{i}$ by $l$ coarse elements, as is shown in Figure \ref{fig:mesh}. We use $\omega_i^{+}$ to denote the oversampled neighborhood enlarged from $\omega_i$. A special case is shown in Figure \ref{fig:mesh}, where we just perform the extension in one direction while size is reduced in another direction to reduce the computation.
\section{Construction of approximation space}\label{sec3}
In this section, we will introduce necessary ingredients for the construction of the mixed CEM-GMsFEM. Specifically, we will present the construction of pressure multiscale basis functions and velocity multiscale basis functions.

On fine grid $\mathcal{T}_h$, we apply the lowest order Raviart-Thomas finite element method (RT0). Let $V_h$ and $Q_h$ be the finite element space for the velocity and pressure functions respectively. More precisely, $V_h$ is spanned by Raviart- Thomas basis functions of lowest order and $Q_h$ is spanned by piecewise constant functions. We want to seek $\mathbf{v}_h\in V_{h}$ and $q_h\in Q_{h}$ as reference solutions by solving the following equations:
 \begin{eqnarray}
 &a(\mathbf{v}_h,\mathbf{w})-b(\mathbf{w},p_h) = 0 \quad\forall \mathbf{w}\in V_h,\label{v_fine}\\
 &(\dot{p}_h,q)_{\rho}+b(\mathbf{v}_h,q)=(f,q)\quad\forall q\in Q_h.\label{q_fine}
 \end{eqnarray}
 The well-poseness of (\ref{v_fine})-(\ref{q_fine}) is proved in \cite{verfurth1984error}.

Then we introduce the mixed CEM-GMsFEM. The key idea in multiscale construction is that we solve local basis functions in some local regions, when the fine-scale information of permeability field is taken into consideration. Due to the localization property, one can obtain global multiscale basis functions based on the local version.

\subsection{Construction of pressure multiscale functions}
We first introduce some notations for later use. Let $S$ be a domain and $V_{h,0}(S)=\{\mathbf{v}\in V_h\cap H(\text{div},S)|\mathbf{v}\cdot \mathbf{n}=0\;\mbox{on}\;\partial S\}$, $Q_{h}(S)=Q_{h}\cap L^2(S).$

Under this setting, we solve a local spectral problem. For each $K_i$ with $i=1,\ldots,N_e$, we seek $(\bm{\phi}_j^i,p_j^i)\in V_{h,0}(K_i)\times Q_{h}(K_i)$ such that
\begin{eqnarray}
\begin{aligned}
a(\bm{\phi}_j^i,\mathbf{v})-b(\mathbf{v},p_j^i)&=0\quad\forall \mathbf{v}\in V_{h,0}(K_i),\\
b(\bm{\phi}_j^i,q)&=\lambda_j^i s_i(p_j^i,q)\quad \forall q \in Q_{h}(K_i),
\label{local spectral pro}
\end{aligned}
\end{eqnarray}
where
\begin{eqnarray*}
s_{i}(p, q)=\int_{K_{i}} \tilde{\kappa} p q, \quad \tilde{\kappa}=\kappa \sum_{j=1}^{N_{c}}\left|\nabla \chi_{j}\right|^{2}.
\end{eqnarray*}
Here, $\{\chi_i\}$ is a set of standard partition of unity functions. $N_c$ is the total number of coarse nodes and $N_{\text{in}}$ is the number of interior coarse nodes.
 For each interior coarse node $x_i$ with $i=1,\ldots,N_{\text{in}}$, we solve $\chi_i$ in corresponding coarse neighborhood $\omega_i$. In detail, $\chi_i$ is solved in each coarse element via the following equations \cite{hou1997multiscale},
\begin{eqnarray}
\begin{aligned}
-\nabla\cdot(\kappa \nabla\chi_i)&=0 \quad\text{in }K\subset \omega_i,\\
\chi_i &=g_i\text{  on }\partial K\backslash\partial\omega_i, \\
\chi_i &=0 \text{  on }\partial\omega_i,
\end{aligned}
\end{eqnarray}
where $g_i$ are some polynomial and continuous functions. And we take linear functions for simplicity.
Based on normalization, we can assume $s_i(p_j^i,p_j^i)=1,\quad \forall j=1,\ldots L_i$. We arrange the eigenvalues in (\ref{local spectral pro}) in nondecreasing order $0<\lambda_1^i\leq \ldots\lambda_{L_i}^i$, and we take the first $J_i$ corresponding eigenfunctions $p_j^i$ with $J_i\leq L_i$ to reduce dimension.

Let the local pressure multiscale space $Q_{\text{ms,i}}=\text{span}\{p_j^i,j=1,\ldots J_i\}$ and further we can obtain $Q_{\text{ms}}=\oplus_{i}^{N_e} Q_{\text{ms,i}}$, which will be the multiscale space for pressure multiscale functions.
\subsection{Construction of velocity multiscale functions}
First of all, we introduce a interpolation operator $\pi:Q\rightarrow Q_{\text{ms}}$ by
\begin{eqnarray}
\pi(q)=\sum_{i=1}^{N} \sum_{j=1}^{J_{i}} \frac{s_{i}\left(q, p_{j}^{(i)}\right)}{s_{i}\left(p_{j}^{(i)},p_{j}^{(i)}\right)} p_{j}^{(i)}  \quad \forall q \in Q.
\end{eqnarray}
For each global function $p_j^i\in Q_{\text{ms}} $, we can obtain a corresponding velocity $\bm{\psi}_j^i$ supported in an oversampled region $K_i^{+}$. Applying similar procedure in last subsection, we can take advantage of localization property of global basis function $\bm{\psi}_j^i$ and solve corresponding local velocity multiscale function $\bm{\psi}_{j,\text{ms}}^i$ in $K_i^{+}$. In practice, as mentioned in \cite{chung2018constraint}, there are two types of ways of velocity multiscale basis functions. In this paper, we just apply the second kind, which gives better accuracy. In particular, for each $p_j^i$ solved in (\ref{local spectral pro}) we can compute $\bm{\psi}_{j,ms}^{(i)}\in V_0(K_i^+)$ and $q_{j,ms}^{(i)}\in Q(K_i^{+})$ such that
\begin{eqnarray}
\begin{aligned}
a(\bm{\psi}_{j,\text{ms}}^i,\mathbf{v})-b(\mathbf{v},q_{j,ms}^i)&=0\quad \forall \mathbf{v} \in  V_0(K_i^{+}),\\
s(\pi q_{j,ms}^i,\pi q)+b(\bm{\psi}_{j,ms}^i,q)&=s(p_j^i,q)\quad \forall q \in Q(K_i^{+}).
\label{cem}
\end{aligned}
\end{eqnarray}
We remark that $\bm{\psi}_{j,ms}^{(i)}\in V_0(K_i^+)$ is defined in the local region $K_i^{+}$ because of the localization property of the global basis function $\bm{\psi}_{j}^{(i)}\in V_0$. More specifically, $\bm{\psi}_{j}^{(i)}$ is obtained by solving the following problem:
\begin{eqnarray}
\begin{aligned}
a(\bm{\psi}_{j}^i,\mathbf{v})-b(\mathbf{v},q_{j}^i)&=0\quad \forall \mathbf{v} \in  V_0(K_i^{+}),\\
s(\pi q_{j}^i,\pi q)+b(\bm{\psi}_{j}^i,q)&=s(p_j^i,q)\quad \forall q \in Q(K_i^{+}).
\label{cem_glo}
\end{aligned}
\end{eqnarray}
Then, we define the multiscale space for velocity as
$$V_{ms}:=\text{span}\{\bm{\psi}_{j}^i: i=1,\ldots,N_e, j=1,\ldots,J_i\}.$$
Then the mixed CEM-GMsFEM for (\ref{model}) reads: Find $(\mathbf{v}_{\text{ms}},p_{\text{ms}})\in V_{\text{ms}}\times Q_{\text{ms}}$ such that
 \begin{eqnarray}
 &a(\mathbf{v}_{\text{ms}},\mathbf{w})-b(\mathbf{w},p_{\text{ms}}) = 0 \quad\forall \mathbf{w}\in V_{\text{ms}} \label{v_ms},\\
 &(\dot{p}_{\text{ms}},q)_{\rho}+b(\mathbf{v}_{\text{ms}},q)=(f,q)\quad\forall q\in Q_{\text{ms}}.\label{q_ms}
 \end{eqnarray}
 The well-poseness of solution to (\ref{v_ms})-(\ref{q_ms}) is proved in Lemma 8 of \cite{chung2018constraint}.
\subsection{Matrix formulation}
To help the readers better understand the proposed method, we will present the following discretization of the concerned problem.
Assuming we obtain the multiscale spaces both for velocity and pressure, in particular,
\begin{eqnarray*}
V_{\text{ms}}=\span\{\bm{\psi}_i\}_{i=1}^S \quad Q_{\text{ms}}=\span\{p_i\}_{i=1}^M,
\end{eqnarray*}
where $M=\sum_{i=1}^{N_e} J_i$, is the dimension for each multiscale space. An equivalent matrix formulation of
(\ref{v_ms})-(\ref{q_ms}) is written as follows:
\begin{eqnarray}
\begin{aligned}
\mathcal{A}\mathbf{v}-\mathcal{B}\mathbf{p}&=0,\\
\mathbf{\mathcal{M}\dot{p}}+\mathcal{B}^{T}\mathbf{v}&=\mathbf{f},
\label{Matrix}
\end{aligned}
\end{eqnarray}
where we have the following definitions of matrices:
\begin{eqnarray*}
\begin{array}{c}
\mathcal{A}:=\left(a\left(\bm{\psi}_{i}, \bm{\psi}_{j}\right)\right) \in \mathbb{R}^{S \times S}, \quad \mathcal{M}:=\left(\left(p_{i}, p_{j}\right)\right) \in \mathbb{R}^{M \times M}, \\
\mathcal{B}:=\left(b\left(\bm{\psi}_{i}, p_{j}\right)\right) \in \mathbb{R}^{S \times M}, \quad \text { and } \quad \mathbf{f}:=\left(\left(f, p_{i}\right)\right) \in \mathbb{R}^{M}.
\end{array}
\end{eqnarray*}
As usual procedure, we use linear combination of basis functions to express solution. $\mathbf{v}:=v(t)=(v_i(t))_{i=1}^{S}$, and $\mathbf{p}:=p(t)=(p_i(t))_{i=1}^{M}$ are vectors of coefficients for the approximations $v_{\text{ms}}$ and $p_{\text{ms}}$. More precisely,
\begin{eqnarray*}
v_{\text{ms}}=\sum_{i=1}^{S}v_i(t)\bm{\psi}_i,\quad
p_{\text{ms}}=\sum_{i=1}^{M}p_i(t)p_i.
\end{eqnarray*}
To discretize in time, we use the backward Euler method, in particular, (\ref{Matrix}) is transformed into
\begin{eqnarray*}
\begin{aligned}
\mathcal{A}\mathbf{v}^{(n+1)}-\mathcal{B}\mathbf{p}^{(n+1)}&=0,\\
\mathbf{\mathcal{M}}\dfrac{\mathbf{p}^{(n+1)}-\mathbf{p}^{(n)}}{\mathbf{\tau}}
+\mathcal{B}^{T}\mathbf{v}^{(n+1)}&=\mathbf{f}^{(n+1)},
\end{aligned}
\end{eqnarray*}
where we have $\mathbf{v}^{(n)}=\mathbf{v}(t_n)$, $\mathbf{p}^{(n)}=\mathbf{p}(t_n)$ and $\mathbf{f}^{(n)}=\left(f(t_n),p_i\right)_{i=1}^{M}$. Here $\tau$ is the time stepsize.

\section{Stability and convergence analysis}\label{analysis}
 In this section, we will prove the stability and derive an error bound for our approximation.
 We use the pair $(\mathbf{v}_h,p_h)$ solved in (\ref{v_fine})-(\ref{q_fine}) as a pair of reference solution. Besides, we define the contrast of the permeability field to be $R_{\kappa}=\frac{\max_{x\in \Omega}\kappa(x)}{\min_{x\in \Omega}\kappa(x)}$ and $J$ to be the number of oversampling layers we enlarge in each oversampled region.

 \subsection{Stability}
 In this section, we prove the stability of the multiscale solution. First of all, we define the norms $\|\cdot\|_a$ for $V_0$ and $\|\cdot\|_{\rho}$ for space $Q$ by
 \begin{equation*}
     \|\mathbf{u}\|_{a}=\int_{\Omega}\kappa^{-1}\mathbf{u}^2 ,\quad \|p\|_{\rho}=\int_{\Omega}\rho p^2,\quad \|p\|_{s}=\int_{\Omega}\tilde{\kappa} p^2.
 \end{equation*}
 \begin{theorem}
 \label{stat}
 Let $(\mathbf{v}_{\text{ms}},p_{\text{ms}})\in V_{\text{ms}}\times Q_{\text{ms}}$ be the solution of (\ref{v_ms})-(\ref{q_ms}). Then we have the following property:
 \begin{eqnarray}
 2\|\mathbf{v}_{\text{ms}}\|_a^2+\dfrac{d}{dt}\|p_{\text{ms}}\|_{\rho}^2=0 \text{ if }f\equiv 0.\label{conser1}
 \end{eqnarray}
  Moreover, the following estimate holds:
 \begin{eqnarray*}
 \begin{aligned}
 \max _{0 \leq t \leq T}\left\|p_{\operatorname{ms}}(t, \cdot)\right\|_{\rho}^{2} &\leq 4\left(\left\|h_{p}\right\|_{\rho}^{2}+\int_{0}^{T}\left\|\rho^{-1} f\right\|_{\rho}^{2} d t\right),\\
 \max _{0 \leq t \leq T}\left\|\mathbf{v}_{\operatorname{ms}}(t, \cdot)\right\|_{a}^{2}&\leq
 2\left(\left\|h_{p}\right\|_{\rho}+\int_{0}^{T}\left\|\rho^{-1} f\right\|_{\rho} d t\right)*\max _{0 \leq t \leq T}\left\|\rho^{-1} f\right\|_{\rho}.
 \end{aligned}
 \end{eqnarray*}
 \end{theorem}
 \begin{proof}
  Taking $\mathbf{w}=\mathbf{v}_{\text{ms}}$ and $p=p_{\text{ms}}$ in (\ref{v_ms})-(\ref{q_ms}) leads to
 \begin{eqnarray}
 2\|\mathbf{v}_{\text{ms}}\|_a^2+\dfrac{d}{dt}\|p_{\text{ms}}\|_{\rho}^2=\left(\rho^{-1}f,p_{\text{ms}}\right)_{\rho},\label{conser2}
 \end{eqnarray}
 which implies (\ref{conser1}).
 Integrating with respect to $t$ and we obtain:
\begin{eqnarray*}
\left\|p_{\mathrm{ms}}(t, \cdot)\right\|_{\rho}^{2} \leq 2 \max _{0 \leq t \leq T}\left\|p_{\mathrm{ms}}(t, \cdot)\right\|_{\rho}\left(\int_{0}^{T}\left\|\rho^{-1} f\right\|_{\rho} d t\right)+\left\|h_{p}\right\|_{\rho}^{2}.
\end{eqnarray*}
Using Young's inequality, we obtain
\begin{eqnarray}
\begin{aligned}
\left\|p_{\operatorname{ms}}(t, \cdot)\right\|_{\rho}^{2} &\leq \frac{1}{2} \max _{0 \leq t \leq T}\left\|p_{\max }(t, \cdot)\right\|_{\rho}^{2}+2\left(\int_{0}^{T}\left\|\rho^{-1} f\right\|_{\rho} d t\right)^{2}+\left\|h_{p}\right\|_{p}^{2} ,\\
\Longrightarrow\max _{0 \leq t \leq T}\left\|p_{\operatorname{ms}}(t, \cdot)\right\|_{\rho}^{2} &\leq 4\left(\left\|h_{p}\right\|_{\rho}^{2}+\int_{0}^{T}\left\|\rho^{-1} f\right\|_{\rho}^{2} d t\right),\\
\max _{0 \leq t \leq T}\left\|p_{\operatorname{ms}}(t, \cdot)\right\|_{\rho}&\leq 2\left(\left\|h_{p}\right\|_{\rho}+\int_{0}^{T}\left\|\rho^{-1} f\right\|_{\rho} d t\right)\label{p_norm}.
\end{aligned}
\end{eqnarray}
Then we can infer from (\ref{conser2}) and (\ref{p_norm}),
\begin{eqnarray}
\begin{aligned}
\max _{0 \leq t \leq T}\left\|\mathbf{v}_{\operatorname{ms}}(t, \cdot)\right\|_{a}^{2}&\leq \max _{0 \leq t \leq T}\left\|\rho^{-1} f\right\|_{\rho}* \max _{0 \leq t \leq T}\left\|p_{\operatorname{ms}}(t, \cdot)\right\|_{\rho}\\
&\leq 2\max _{0 \leq t \leq T}\left\|\rho^{-1} f\right\|_{\rho}*\left(\left\|h_{p}\right\|_{\rho}+\int_{0}^{T}\left\|\rho^{-1} f\right\|_{\rho} d t\right).
\end{aligned}
\end{eqnarray}
This strategy will be applied in the proof of Theorem \ref{thm2}.
 \end{proof}
 \subsection{Error analysis}
 Before we turn to the error analysis section, we first define the multiscale projection.  We consider the pair $(\mathbf{v},p)\in V_h\times Q_h$ satisfying
\begin{eqnarray*}
\begin{aligned}
a(\mathbf{v},\mathbf{w})-b(\mathbf{w},p)&=0\quad\forall \mathbf{w}\in V_{h},\\
b(\mathbf{w},q)&=(\alpha,q) \quad\forall q\in Q_{h},
\end{aligned}
\end{eqnarray*}
for some $\alpha\in L^{2}(\Omega)$. We define $\sigma \mathbf{v}\in V_{\text{ms}}$ and $\sigma p\in Q_{\text{ms}}$ to be the multiscale projection of $(\mathbf{v},p)$ if they satisfy
\begin{eqnarray}
\begin{aligned}
a(\sigma \mathbf{v},\mathbf{w})-b(\mathbf{w},\sigma p)&=0\quad\forall \mathbf{w}\in V_{\text{ms}},\\
b(\sigma \mathbf{v},q)&=(\alpha,q) \quad\forall q\in Q_{\text{ms}}. \label{proj}
\end{aligned}
\end{eqnarray}
 \begin{lemma}
Assume that $J=O(\log(R_{\kappa}/H^2))$. For any $\mathbf{v} \in V_{h},$ the following estimate holds:
\begin{eqnarray}
\|\mathbf{v}-\sigma \mathbf{v}\|_{a} \leq CH\Lambda^{-1 / 2} \label{v_inter},
\end{eqnarray}
where $\Lambda:=\min _{1 \leq i \leq N_e} \lambda_{J_{i}+1}^{i}$ and $\left\{\lambda_{j}^{i}\right\}$ are the eigenvalues obtained from the construction of multiscale basis functions.
 \end{lemma}
\begin{proof}
 For any $\mathbf{v} \in V_{h},$ we define $\beta \in Q_{h}\left(\text { with } \int_{\Omega} \beta d x=0\right)$ such that
\begin{eqnarray*}
b(\mathbf{w}, \beta)=a(\mathbf{v}-\sigma \mathbf{v}, \mathbf{w}) \quad \forall \mathbf{w} \in V_{h},
\end{eqnarray*}
where the existence and uniqueness of such $\beta$ can be guaranteed by
(\ref{inf-sup}).
Denote $\mathbf{z}=\mathbf{v}-\sigma \mathbf{v}$. Then, $(\mathbf{z}, \beta) \in V_{h} \times Q_{h}$ satisfies the following system:
\begin{eqnarray*}
\begin{aligned}
a(\mathbf{z}, \mathbf{w})-b(\mathbf{w}, \beta) &=0 & & \forall \mathbf{w} \in V_{h}, \\
b(\mathbf{z}, q) &=(\nabla \cdot(\mathbf{v}-\sigma \mathbf{v}), q) & & \forall q \in Q_{h}.
\end{aligned}
\end{eqnarray*}

Hence, the following estimate holds:
\begin{eqnarray*}
\|\mathbf{z}-\sigma \mathbf{z}\|_{a} \leq CH\|\nabla \cdot(\mathbf{z}-\sigma \mathbf{z})\|_{L^{2}(\Omega)}.
\label{z_error}
\end{eqnarray*}
using the result of Theorem 1 in \cite{chung2018constraint}. We emphasize here that (\ref{z_error}) is under the assumption $J=O(\log(R_{\kappa}/H^2))$, which will be illustrated in the final remark of this section.  Therefore, we have
\begin{eqnarray*}
\begin{aligned}
\|\mathbf{v}-\sigma \mathbf{v}\|_{a}^{2} &=a(\mathbf{z}, \mathbf{v}-\sigma \mathbf{v})=a(\mathbf{z}-\sigma \mathbf{z}, \mathbf{v}-\sigma \mathbf{v}) \\
& \leq\|\mathbf{z}-\sigma \mathbf{z}\|_{a}\|\mathbf{v}-\sigma \mathbf{v}\|_{a} \\
& \leq C H\|\nabla \cdot(\mathbf{v}-\sigma \mathbf{v})\|_{L^{2}(\Omega)}\|\mathbf{v}-\sigma \mathbf{v}\|_{a},
\end{aligned}
\end{eqnarray*}
where we make use of $\sigma \mathbf{z}=0$.
Based on the construction of the basis construction (one may refer to (\ref{local spectral pro}) and (\ref{cem})), since $b(\mathbf{v}-\sigma \mathbf{v}, q)=s\left(\tilde{\kappa}^{-1} \nabla \cdot(\mathbf{v}-\sigma \mathbf{v}), q\right)=0$ for all $q \in Q_{\mathrm{ms}},$
there exists a set of real numbers $\left\{c_{j}^{i}\right\}$ such that
\begin{eqnarray*}
\tilde{\kappa}^{-1} \nabla \cdot(\mathbf{v}-\sigma \mathbf{v})=\sum_{i=1}^{N} \sum_{j>J_{i}} c_{j}^{i} p_{j}^{i}.
\end{eqnarray*}
Then, by the orthogonality of the eigenfunctions $\left\{p_{j}^{i}\right\}$ and the spectral problem (\ref{local spectral pro}), we have
\begin{eqnarray*}
\begin{aligned}
\|\nabla \cdot(\mathbf{v}-\sigma \mathbf{v})\|_{L^{2}(\Omega)}^{2} &=\left\|\tilde{\kappa}^{-1 / 2} \nabla \cdot(\mathbf{v}-\sigma \mathbf{v})\right\|_{s}^{2} \\
& = \sum_{i=1}^{N} \sum_{j>J_{i}}\left(c_{j}^{i}\right)^{2}\left\|p_{j}^{i}\right\|_{s}^{2} \leq C \Lambda^{-1} \sum_{i=1}^{N} \sum_{j>J_{i}}\left(c_{j}^{i}\right)^{2}\left\|\phi_{j}^{i}\right\|_{a}^{2},
\end{aligned}
\end{eqnarray*}
which proves (\ref{v_inter}).
\end{proof}
\begin{theorem}
Assume that $J=O(\log(R_{\kappa}/H^2))$. Suppose that $(\mathbf{v}_h,p_h)$ is the solution to (\ref{v_fine})-(\ref{q_fine}),$(\sigma \mathbf{v}_h,\sigma p_h)$ is the multiscale projection of the $(\mathbf{v}_h,p_h)$. $(v_{\text{ms}},p_{\text{ms}})$ is the solution to (\ref{v_ms})-(\ref{q_ms}). Then, the following estimate holds:
\begin{eqnarray*}
\|\sigma \mathbf{v}_h- \mathbf{v}_{\text{ms}}\|_a^2\leq C\left(\|\tilde{\kappa}^{-1}F\|_{s}^2+\int_{0}^T\|\tilde{\kappa}^{-1}\Dot{F}\|_{s}^2 dt\right).
\end{eqnarray*}
When $\{\chi_i\}$ is the standard bilinear functions, we have
\begin{eqnarray}
\|\sigma \mathbf{v}_h- \mathbf{v}_{\text{ms}}\|_a^2\leq CH^2(\int_{0}^{T} \|\dot{F}\|_{L^2}^2 dt+ \|\dot{F}\|_{L^2}^2).
\end{eqnarray}
Combing (\ref{v_inter}), we have obtained
\begin{eqnarray}
\| \mathbf{v}_h- \mathbf{v}_{\text{ms}}\|_a^2\leq CH^2(\int_{0}^{T} \|\dot{F}\|_{L^2}^2 dt+ \|\dot{F}\|_{L^2}^2+\frac{1}{\Lambda}).
\end{eqnarray}
\label{thm2}
\end{theorem}
\begin{proof}
Combining the definition of multiscale projection and (\ref{v_fine})-(\ref{q_fine}), we have
 \begin{eqnarray}
 \begin{aligned}
   a(\mathbf{v}_h-\sigma \mathbf{v}_h,\mathbf{w})-b(\mathbf{w},p_h-\sigma p_h) &= 0 \quad\forall \mathbf{w}\in V_h,\\
b(\mathbf{v}_h-\sigma \mathbf{v}_h,q)&=(F,q)\quad\forall q\in Q_h, \label{diff1}
 \end{aligned}
 \end{eqnarray}
 where $F=f-\rho p_h$. By the result of Theorem 1 in \cite{chung2018constraint}, we have the following error estimate.
 \begin{eqnarray}
     \|\mathbf{v}_h-\sigma \mathbf{v}_h\|_a^2 + \|p_h-\sigma p_h\|_{\rho}^2\leq C\|\tilde{\kappa}^{-1}F\|_{s}^2 \label{estimate1}.
 \end{eqnarray}
 Differentiate both two equations in (\ref{diff1}) w.r.t $t$, we have
  \begin{eqnarray*}
 \begin{aligned}
   a(\dot{\mathbf{v}_h}-\dot{\sigma \mathbf{v}_h},\mathbf{w})-b(\mathbf{w},\dot{p_h}-\dot{\sigma p_h}) &= 0 \quad\forall \mathbf{w}\in V_h,\\
b(\dot{\mathbf{v}_h}-\dot{\sigma \mathbf{v}_h},q)&=(\dot{F},q)\quad\forall q\in Q_h.
 \end{aligned}
 \end{eqnarray*}
 Similarly, we can get
 \begin{eqnarray}
      \|\dot{\mathbf{v}_h}-\dot{\sigma \mathbf{v}_h}\|_a^2 + \|\dot{p_h}-\dot{\sigma p_h}\|_{\rho}^2 \leq C\|\tilde{\kappa}^{-1}\Dot{F}\|_{s}^2.
      \label{estimate2}
 \end{eqnarray}
Based on (\ref{v_fine})-(\ref{q_ms}), we have
 \begin{eqnarray*}
 \begin{aligned}
   a(\mathbf{v}_{\text{ms}}-\mathbf{v}_h,\mathbf{w})-b(\mathbf{w},p_{\text{ms}}-p_h) &= 0 \quad\forall \mathbf{w}\in V_h, \\
 (\dot{p_{\text{ms}}}-\dot{p_h},q)_{\rho}+b(\mathbf{v}_{\text{ms}}-\mathbf{v}_h,q)&=0\quad\forall q\in Q_h,
 \end{aligned}
 \end{eqnarray*}
which is equivalent to
 \begin{eqnarray}
 \begin{aligned}
   a(\mathbf{v}_{\text{ms}}-\sigma \mathbf{v}_h,\mathbf{w})-b(\mathbf{w},p_{\text{ms}}-\sigma p_h) &= a(\mathbf{v}_h-\sigma \mathbf{v}_h,\mathbf{w})-b(\mathbf{w},p_{h}-\sigma p_h) \quad\forall \mathbf{w}\in V_h, \\
(\dot{p_{\text{ms}}}-\dot{\sigma p_h},q)_{\rho}+ b(\mathbf{v}_{\text{ms}}-\sigma \mathbf{v}_h,q)&=(\dot{p_{h}}-\dot{\sigma p_h},q)_{\rho}+b(\mathbf{v}_{h}-\sigma \mathbf{v}_h,q)\quad\forall q\in Q_h.\label{add}
 \end{aligned}
 \end{eqnarray}
 Take $\mathbf{w}=\mathbf{v}_{\text{ms}}-\sigma \mathbf{v}_h$ and $q=p_{\text{ms}}-\sigma p_h$, we have
 \begin{eqnarray*}
 LHS=a(\sigma \mathbf{v}_h-\mathbf{v}_{\text{ms}},\sigma \mathbf{v}_h-\mathbf{v}_{\text{ms}})+(\dot{\sigma p_h}-\dot{p_{\text{ms}}},\sigma p_h-p_{\text{ms}})_{\rho}
 \end{eqnarray*}
 and
 \begin{eqnarray*}
    RHS&=a(\sigma \mathbf{v}_h-\mathbf{v}_h,\sigma \mathbf{v}_h-\mathbf{v}_{\text{ms}})-b(\sigma \mathbf{v}_h-\mathbf{v}_{\text{ms}},\sigma p_h-p_h)\\&+(\dot{\sigma p_h}-\dot{p_{h}},\sigma p_h-p_{\text{ms}})_{\rho}+b(\sigma \mathbf{v}_h-\mathbf{v}_h,\sigma p_h-p_{\text{ms}}).
 \end{eqnarray*}
 Using (\ref{proj}), we have
 \begin{eqnarray*}
 \begin{aligned}
 a(\sigma \mathbf{v}_h-\mathbf{v}_h,\sigma \mathbf{v}_h-\mathbf{v}_{\text{ms}})-b(\sigma \mathbf{v}_h-\mathbf{v}_{\text{ms}},\sigma p_h-p_h)&=0,\\
 b(\sigma \mathbf{v}_h-\mathbf{v}_h,\sigma p_h-p_{\text{ms}})&=0.
 \end{aligned}
 \end{eqnarray*}
 Therefore, (\ref{add}) is reduced to
 \begin{eqnarray}
 \|\sigma \mathbf{v}_h-\mathbf{v}_{\text{ms}}\|_a^2+\dfrac{d}{2 dt}\|p_h-p_{\text{ms}}\|_{\rho}^2= (\dot{\sigma p_h}-\dot{p_{h}},\sigma p_h-p_{\text{ms}})_{\rho} \label{add2}.
 \end{eqnarray}
 First consider
  \begin{eqnarray*}
 \dfrac{d}{2 dt}\|p_h-p_{\text{ms}}\|_{\rho}^2\leq (\dot{\sigma p_h}-\dot{p_{h}},\sigma p_h-p_{\text{ms}})_{\rho}.
 \end{eqnarray*}
 Consider the strategy used in Theorem \ref{stat}, we can obtain
 \begin{eqnarray*}
 \max _{0 \leq t \leq T} \|p_h-p_{\text{ms}}\|_{\rho}^2\leq 2 \|p_h|_{t=0}-p_{\text{ms}}|_{t=0}\|_{\rho}^2
 +4 \int_{0}^{T}\|\dot{\sigma p_h}-\dot{p_h}\|_{\rho}^2.
 \end{eqnarray*}
Since we can take $p_h(0,x)=p_{\text{ms}}(0,x)$, we have
 \begin{eqnarray*}
 \max _{0 \leq t \leq T} \|\sigma p_h-p_{\text{ms}}\|_{\rho}^2\leq 4 \int_{0}^{T}\|\dot{\sigma p_h}-\dot{p_h}\|_{\rho}^2.
 \end{eqnarray*}
 By (\ref{estimate2}), we obtain
 \begin{eqnarray*}
 \max _{0 \leq t \leq T} \|\sigma p_h-p_{\text{ms}}\|_{\rho}^2\leq 4  \int_{0}^{T}\|\tilde{\kappa}^{-1}\Dot{F}\|_{s}^2.
 \end{eqnarray*}
 Apply (\ref{add2}) again and (\ref{estimate1}), we have
 \begin{eqnarray*}
   \|\sigma \mathbf{v}_h-\mathbf{v}_{\text{ms}}\|_a^2 \leq 2 \|\dot{\sigma p_h}-\dot{p_{h}}\|_{\rho}^2+\frac{1}{2}
   \|\sigma p_h-p_{\text{ms}}\|_{\rho}^2\leq C\left(\|\tilde{\kappa}^{-1}F\|_{s}^2+\int_{0}^T\|\tilde{\kappa}^{-1}\Dot{F}\|_{s}^2 dt\right).
 \end{eqnarray*}

\end{proof}
We remark that to obtain the $ O(H)$ convergence, we need to choose the size of the oversampling domain $J=O(log(R_{\kappa}/H^2))$, where $H$ is the coarse-mesh size and $J$ is the number of coarse-grid layer extension (see Figure \ref{fig:mesh}). Note that $\|\tilde{\kappa}^{-1}F\|_{s}^2=CH^2/\kappa\leq CH^2 C_{\text{min}}$, where $C_{\text{min}}=\max_{x\in \Omega}\tilde{\kappa}^{-1}(x)$. 
We will discuss the more precise choice of the parameter $J$. In particular, we need to choose the size of the oversampling domain such that (cf. \cite{chung2018constraint})
$$ C_L(J+1)^dE(1+C_{\text{glo}}^2)(1+C_{\text{ms}})^4C_{\text{min}} =O(1),$$
where the constants are defined in (\cite{chung2018constraint}) and are summarized as follows. We have
  $C_L$ is the maximal number of overlapping oversampled regions, $C_{\text{glo}}=O(\max(\tilde{\kappa}))$,
$$C_{\text{ms}}^2=2\dfrac{C_{\text{glo}^2}+\alpha}{(1-\sqrt{\alpha})^2},\quad\displaystyle E=C(1+\frac{1}{\Lambda})(1+C^{-1}(1+\frac{1}{\Lambda})^{-\frac{1}{2}})^{1-J},$$
and $\alpha$ is chosen such that
$$C_L(J+1)^dE(1+C_{\text{glo}}^2)^2\leq \alpha<1.$$
Note that, the constant $E$ decays exponentially with respect to $J$.
For more detail about the derivation of these constants, see (\cite{chung2018constraint}).

\section{Numerical results}\label{sect_numerical}

\begin{figure}[!htbp]
\centering
\subfigure[$\log\kappa_1(x)$]
{ \includegraphics[width=0.45\textwidth]{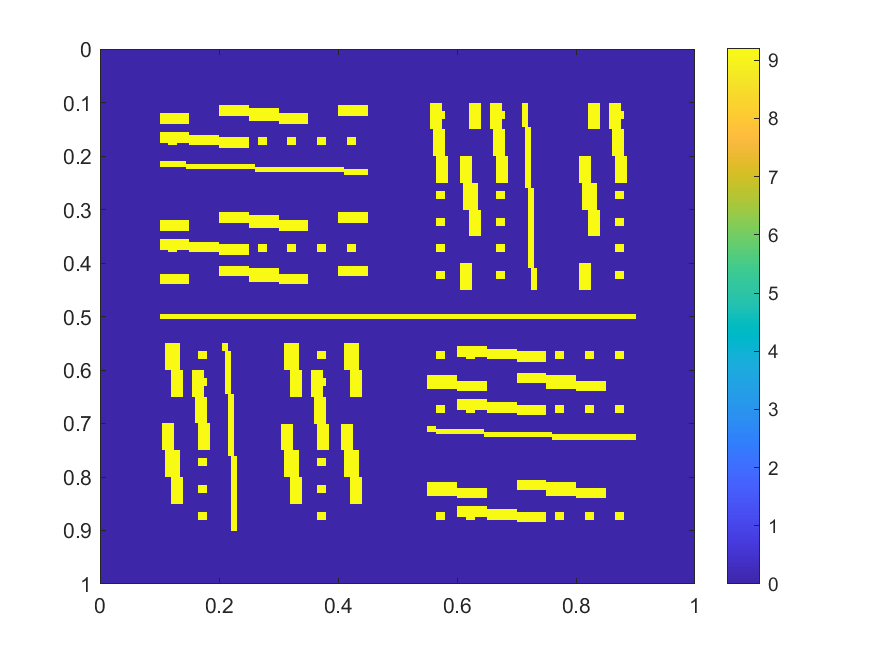}}
\subfigure[$\log\kappa_2(x)$]
{ \includegraphics[width=0.45\textwidth]{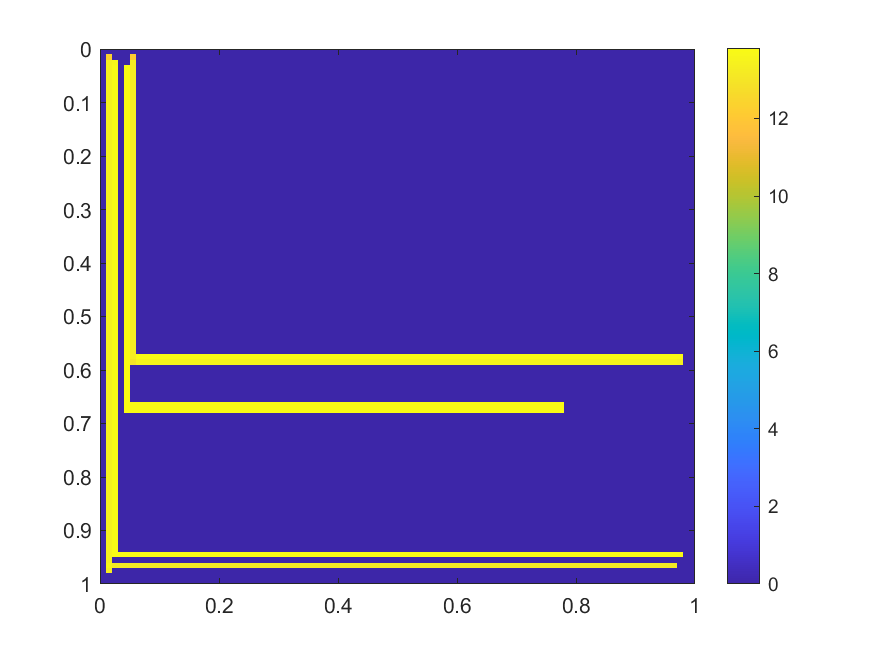}}
\subfigure[$\log\kappa_3(x)$]
{ \includegraphics[width=0.45\textwidth]{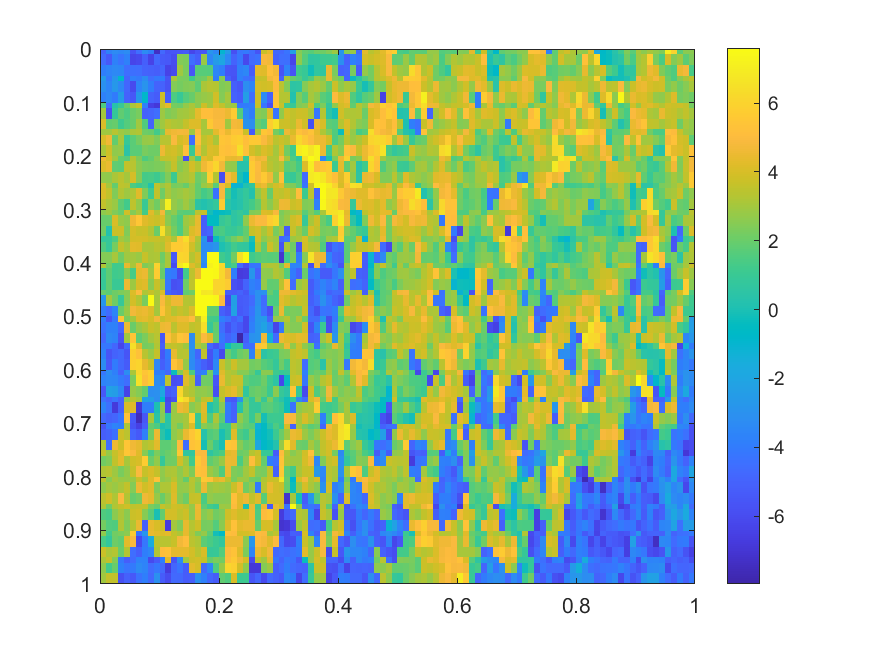}}
\subfigure[$\log\kappa_4(x)$]
{ \includegraphics[width=0.45\textwidth]{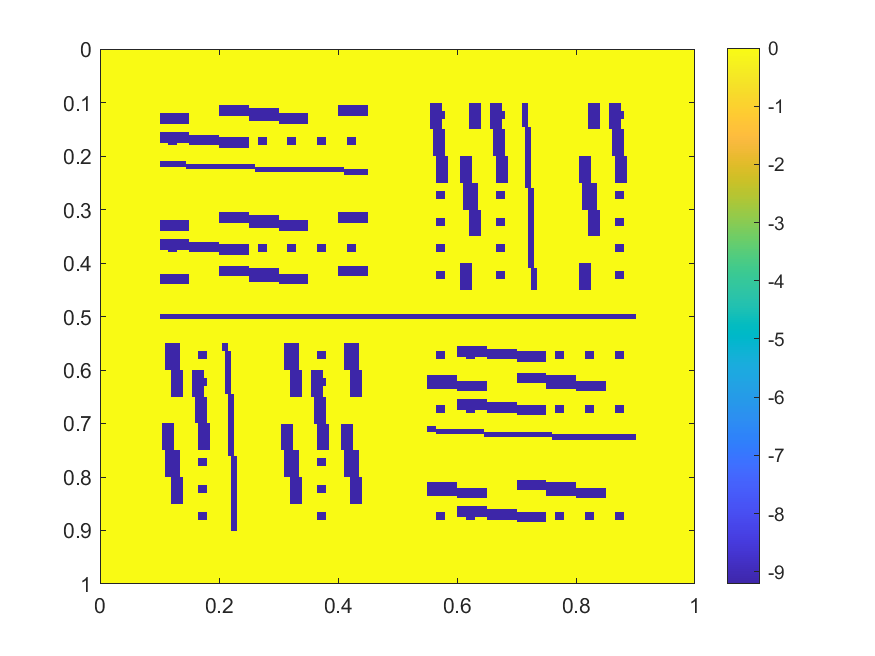}}
\caption{The heterogeneous permeability fields used in our simulations in log scale.}\label{medium}
\end{figure}

In this section, we will present the results of numerical experiments. In this paper, we set $\Omega=[0,1]\times[0,1]$. In addition, we let $h_v=0$ and $h_p=0$. For the source term $f$, we set $f=1$ in $[0,0.1]\times[0,0.1]$ and $f=-1$ in $[0.9,1]\times [0.9,1]$ while 0 elsewhere. We consider four different kinds of permeability fields: $\kappa_1$, $\kappa_2$, $\kappa_3$ and $\kappa_4$, which are shown in Figure \ref{medium}. Note that $\kappa_1$ is a deterministic permeability field with relatively high contrast $10^4$.
It is apparent that abrupt transitions occur between low and high permeability. 
In terms of $\kappa_2$, it is a simplified model for a fractured porous media, which is characterized by complex fracture distribution and high contrast. As for $\kappa_3$, it is extracted from benchmark permeability field: the tenth SPE comparative solution project (SPE10), which is commonly used for assessing upscaling and multiscale methods. Despite no clear separate channels can be detected, $\kappa_3$ is highly heterogeneous and of high contrast which is approximately $10^7$, which creates significant difficulty in computation.  The last one $\kappa_4$ is similar to $\kappa_1$ in view of high contrast and abrupt transitions. However, $\kappa_{\text{max}}=1$ and $0\leq\kappa_{\text{min}}\ll 1$.

For notation, we set $L_z$ to be the number of multiscale (pressure or velocity) basis functions in each coarse element, $J$ is the number of oversampling layers used in (\ref{cem}) and (\ref{cem_glo}) and $H$ is the coarse-mesh size. We remark that each multiscale pressure basis function corresponds to exactly one multiscale velocity basis function.  For each permeability field, we consider different combinations of $L_z$, $J$ and $H$. To illustrate the performance of the method, we use the following error quantities (pressure error: $e_p$, velocity error: $e_v$)
 \begin{equation*}
     e_p=\dfrac{\|p_h-p_{\text{ms}}\|_{L^2(\Omega)}}{\|p_h\|_{L^2(\Omega)}},\quad
     e_v=\dfrac{\|\mathbf{v}_h-\mathbf{v}_{\text{ms}}\|_a}{\|\mathbf{v}_h\|_a},
 \end{equation*}
where $(\mathbf{v}_h,p_h)$ is the reference solution obtained by (\ref{v_fine})-(\ref{q_fine}) on the fine-grid $\mathcal{T}^h$ using the standard mixed finite element method while $(\mathbf{v}_{\text{ms}},p_{\text{ms}})$ is multiscale solution obtained by (\ref{v_ms})-(\ref{q_ms}).

\subsection{Example 1}
\begin{figure}[!htbp]
    \includegraphics[width=6in]{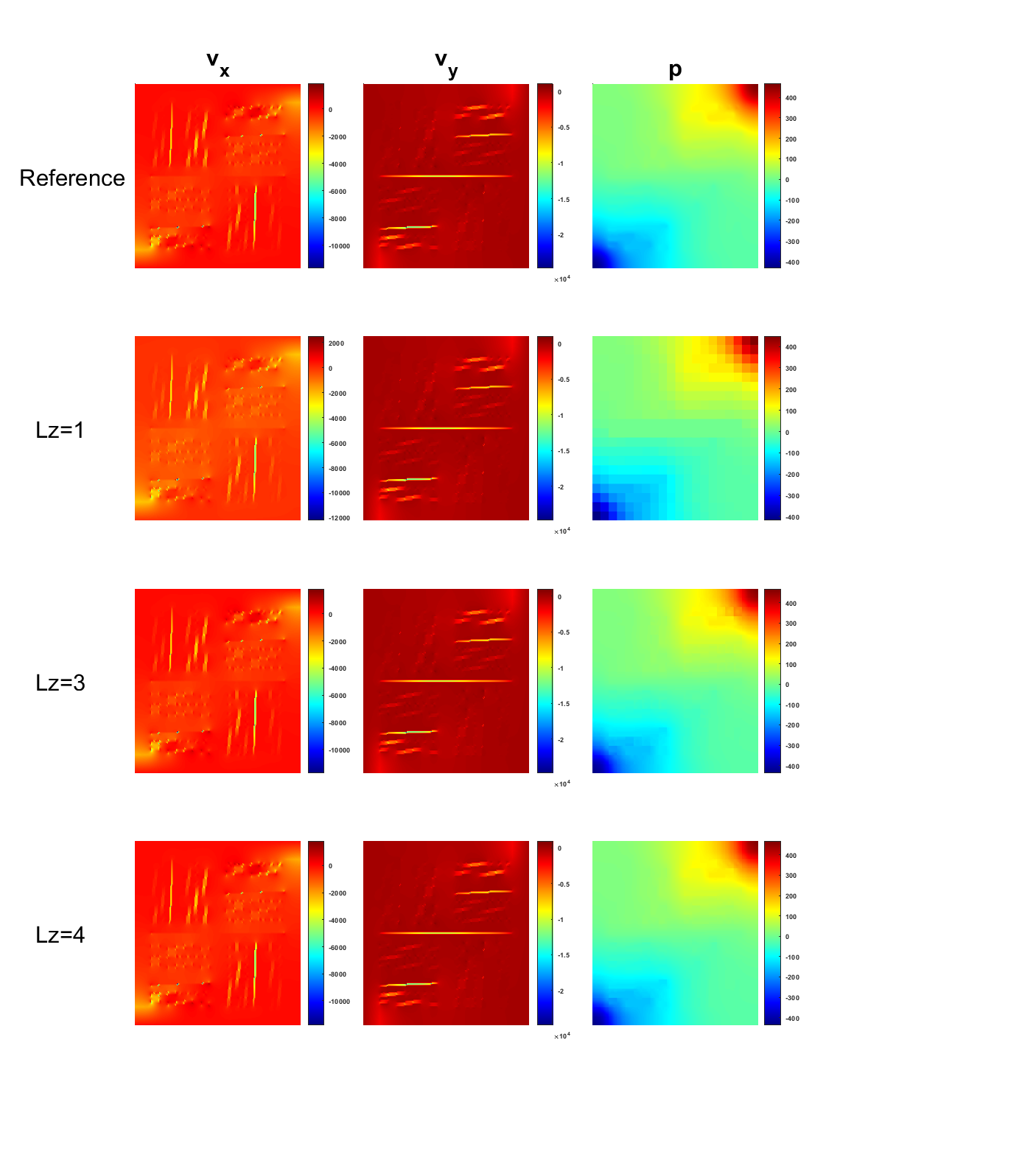}\caption{Example 1. Profile of velocity and pressure with $\mathbf{\kappa_1}$ at $T$ obtained by four methods. $J=2$ and $H=1/20$. Components of velocity solutions in x and y direction are shown in first two columns while pressure solution is in last column. Results from MFEM (reference),  mixed CEM-GMsFEM  with $L_z=1$, $L_z=3$ and $L_z=4$ are shown in one to fourth rows. }
    \label{vp_k1_t3}
\end{figure}

In this part, we apply $\kappa_1$ in (\ref{model}), whose meshsize is $200\times 200$. Besides, we set terminal time $T=1$ and time step size $\tau=10^{-2}$. We compare MFEM (reference method) and mixed CEM-GMsFEM with different settings. In Figure \ref{vp_k1_t3}, we compare reference solution with three cases of mixed CEM-GMsFEM ($L_z=1,3,4$) at terminal time $T$. Each row is a set of solution solved with same setting. In particular, $v_x$ and $v_y$ are the horizontal and vertical components for velocity solution respectively and $p$ is pressure solution. As one can observe from the figure, more basis functions can better approximate the reference case. There are some noticeable mismatches between first two rows. For example, the ranges of $v_x$ are evidently different since the maximum the of second row is much greater than the first row. Besides, some blocks in pressure solution around the lower-left and upper-right corners are easily observed in the second row, which are rarely seen in the last two rows. The last two approximations especially the $L_z=4$ case are almost identical to the reference. It is worth mentioning that even if the multiscale solution is computed in coarse scale, most of the fine-scale details are retained in the approximation with well-chosen $L_z$, which exhibits the power of multiscale methods.

\begin{figure}[!htbp]
\subfigure[velocity error with $L_z$]
{\includegraphics[width=0.48\textwidth]{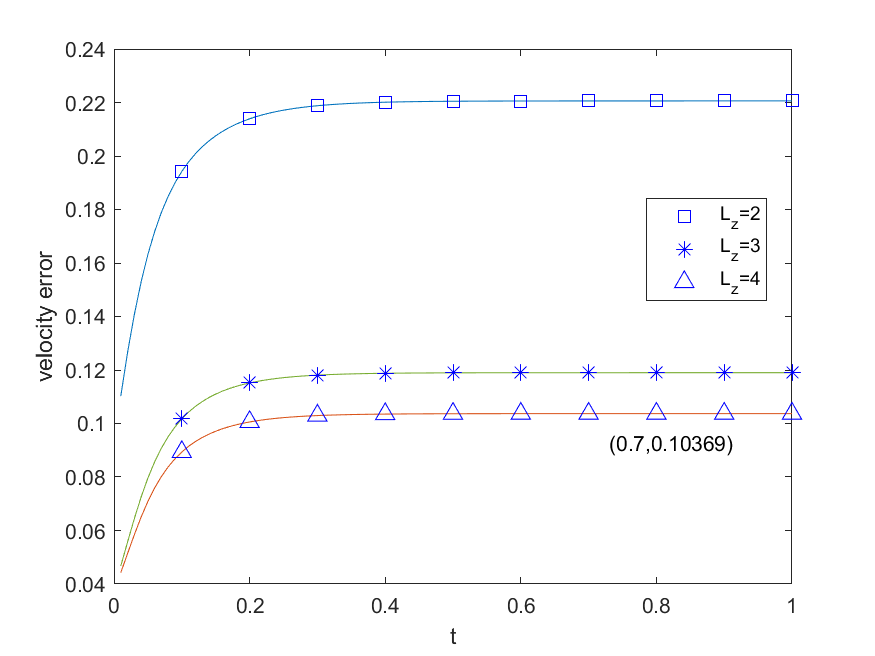}}
\subfigure[pressure error with $L_z$]
{\includegraphics[width=0.48\textwidth]{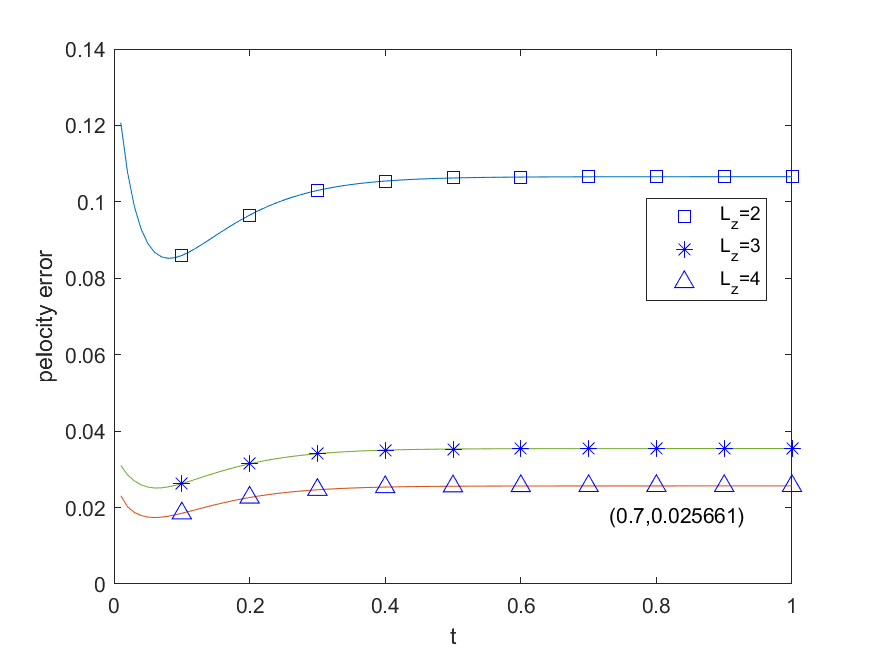}}
\subfigure[velocity error with $J$]
{\includegraphics[width=0.48\textwidth]{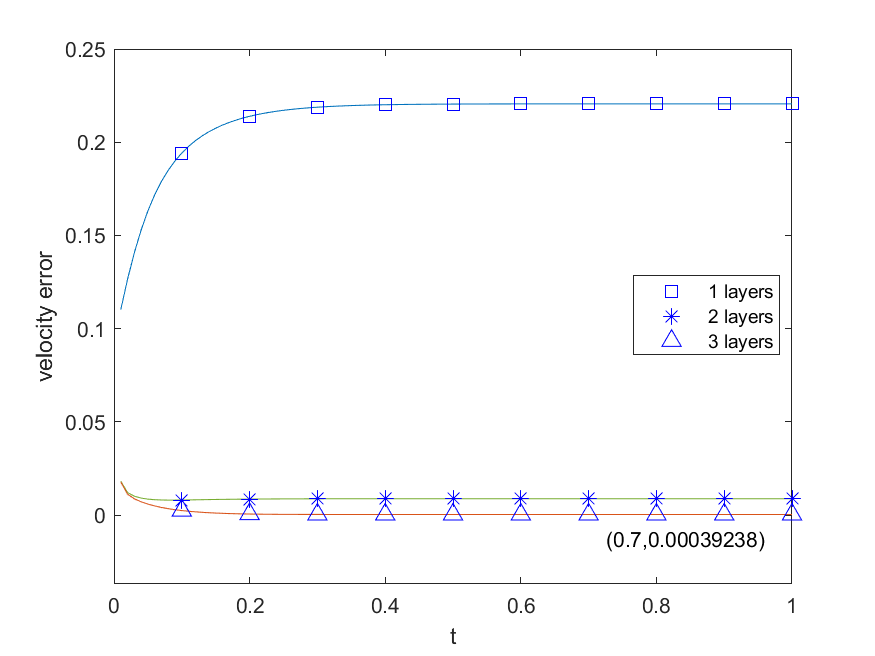}}
\subfigure[pressure error with $J$]
{\includegraphics[width=0.48\textwidth]{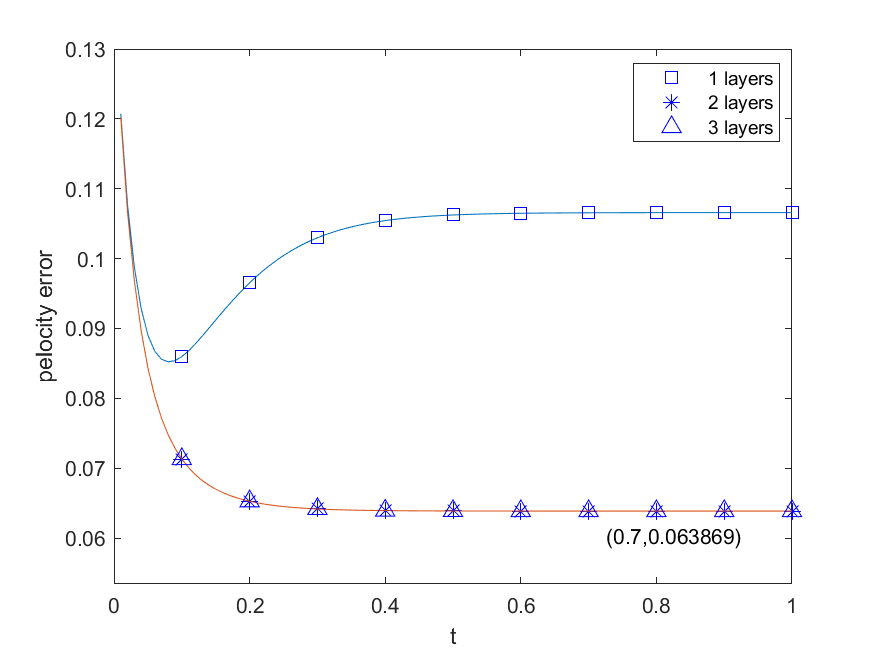}}
\subfigure[velocity error with $H$]
{\includegraphics[width=0.48\textwidth]{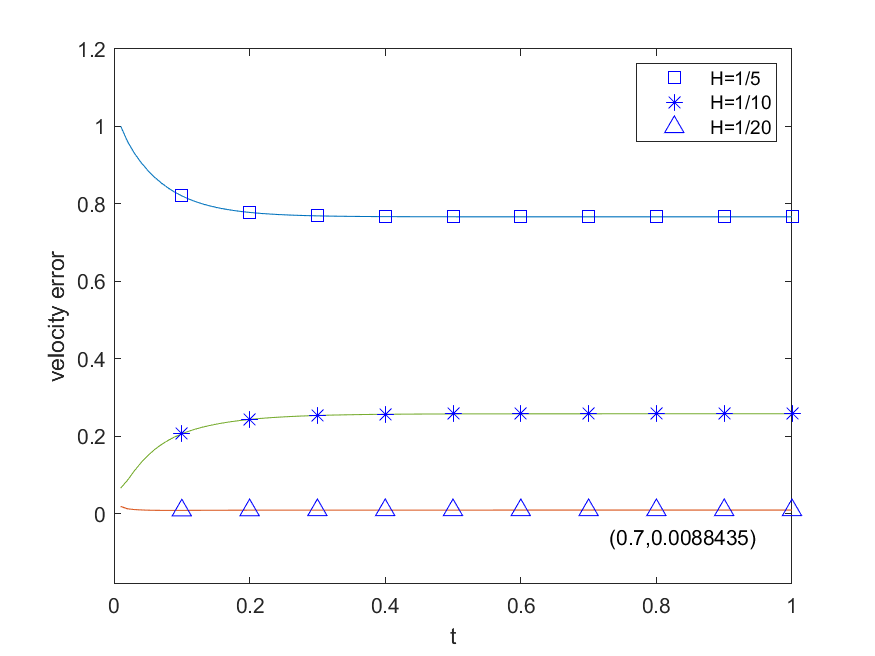}}
\subfigure[pressure error with $H$]
{\includegraphics[width=0.48\textwidth]{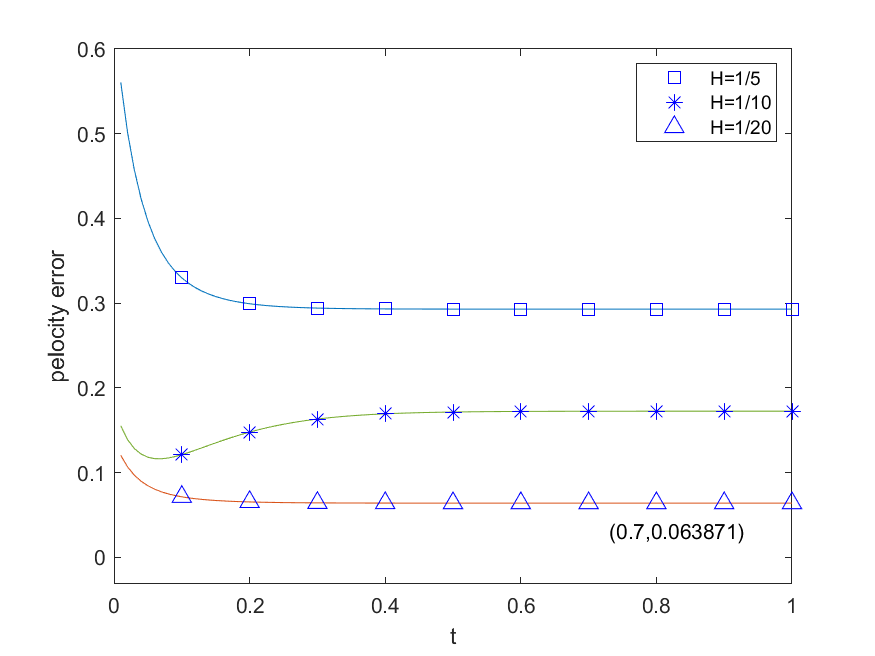}}
\caption{Example 1. Velocity errors and pressure errors with $\mathbf{\kappa_1}$ under different $L_z$, $J$ and $H$. First row : $J=1$, $H=1/20$. Second row: $L_z=2$, $H=1/20$. Third row: $L_z=2$, $J=2$.}\label{error_k1}
\end{figure}

In Figure \ref{error_k1}, we present changes of velocity error and pressure error from initial time to terminal time $T$. In this paper, we mainly consider parameters: $L_z$, $J$ and $H$, the influences of which are shown in three rows accordingly. In each row, we fix two parameters and compare different pairs of velocity and pressure errors by changing the third parameter. In other words, the velocity errors and pressure errors in each row are computed under same setting. In each subfigure, for better representation, we also choose a specific intermediate time step to show the corresponding lowest error at that step. As in (b) in Figure \ref{error_k1}, we choose $t=0.7$, and show the lowest pressure error is about $0.025$ obtained with $L_z=4$. From the first row, we can see the improvements of accuracy by enriching the offline multiscale space from $L_z=2$ to $L_z=4$, where the velocity error is decreased by half from about $22\%$ to $10\%$ and pressure error is largely diminished from $10\%$ to $2\%$. This is consistent to the results in Figure \ref{vp_k1_t3}. In terms of oversampling strategies, one may easily observe from the second row that the accuracy is significantly increased from adding only one layer from $J=1$. Here we choose $L_z=2$ and $H=1/20$. Both of velocity and pressure error are greatly reduced, especially for velocity error, from $20\%$ to $0.04\%$, which shows the efficiency of oversampling strategies. As to the choices of $H$, we consider three cases $1/5$, $1/10$ and $1/20$. One may see the convergence as we decrease $H$. In particular, from $1/5$ to $1/20$, the velocity error is diminished to a great extent from about $80\%$ to less than $1\%$, which is consistent to the convergence analysis in section 4. Due to the independent property of local multiscale space construction, one can construct in parallel to reduce the computation time. Compared to velocity solution, there is a degradation in the accuracy increase in pressure solution but error reduction is still remarkable, which goes from $30\%$ to $6\%$.

\subsection{Example 2}
\begin{figure}[!htbp]
    \includegraphics[width=6in]{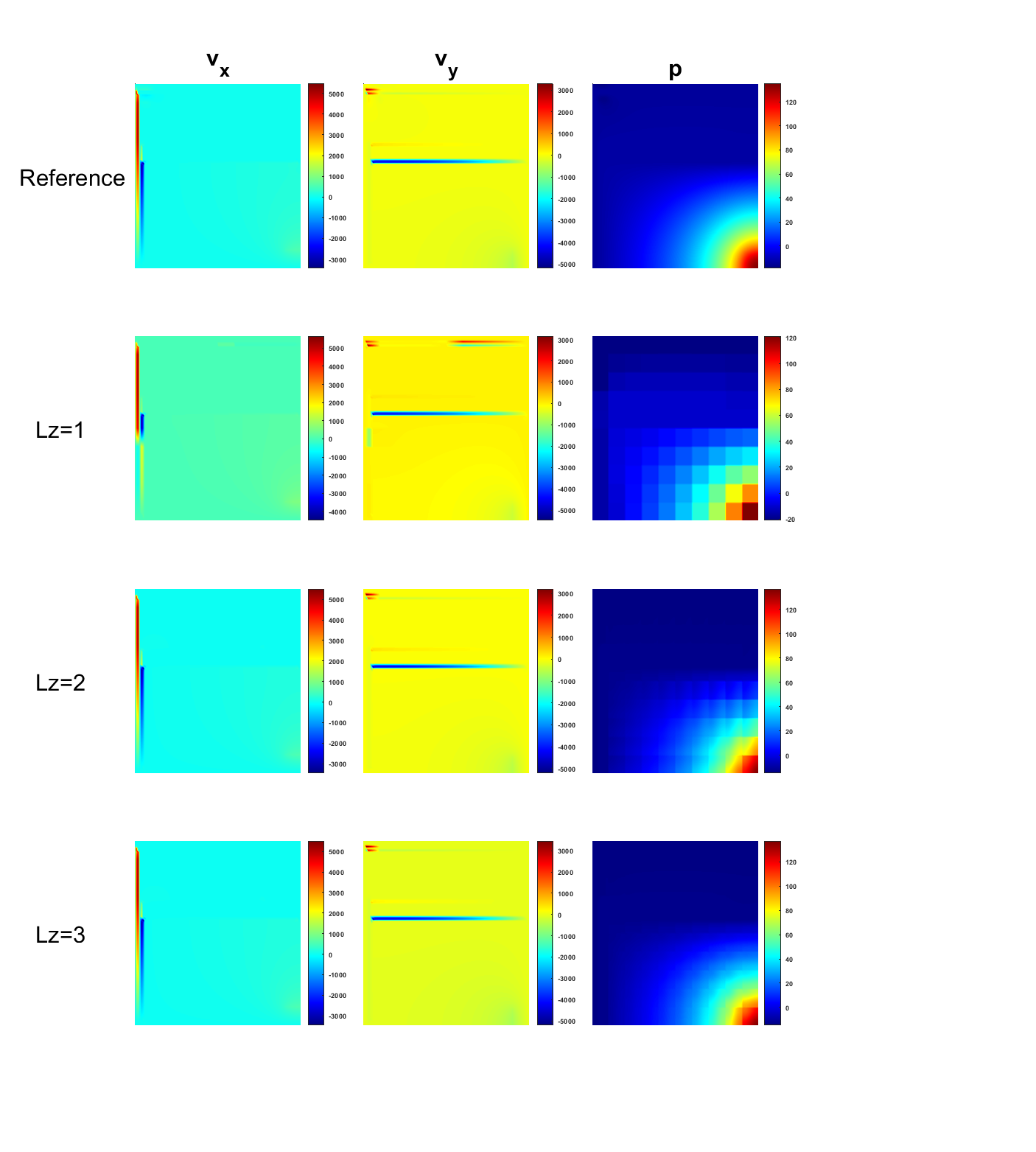}\caption{Example 2. Profile of velocity and pressure with $\mathbf{\kappa_2}$ at $T$ obtained by four methods. $J=3$ and $H=1/10$. Components of velocity solutions in x and y direction are shown in first two columns while pressure solution is in last column. Results from MFEM (reference), mixed CEM-GMsFEM  with $L_z=1$, $L_z=2$ and $L_z=3$ are shown in one to fourth rows.  }
    \label{vp_k2_t3}
\end{figure}

\begin{figure}[!htbp]
\subfigure[velocity error with $L_z$]
{\includegraphics[width=0.48\textwidth]{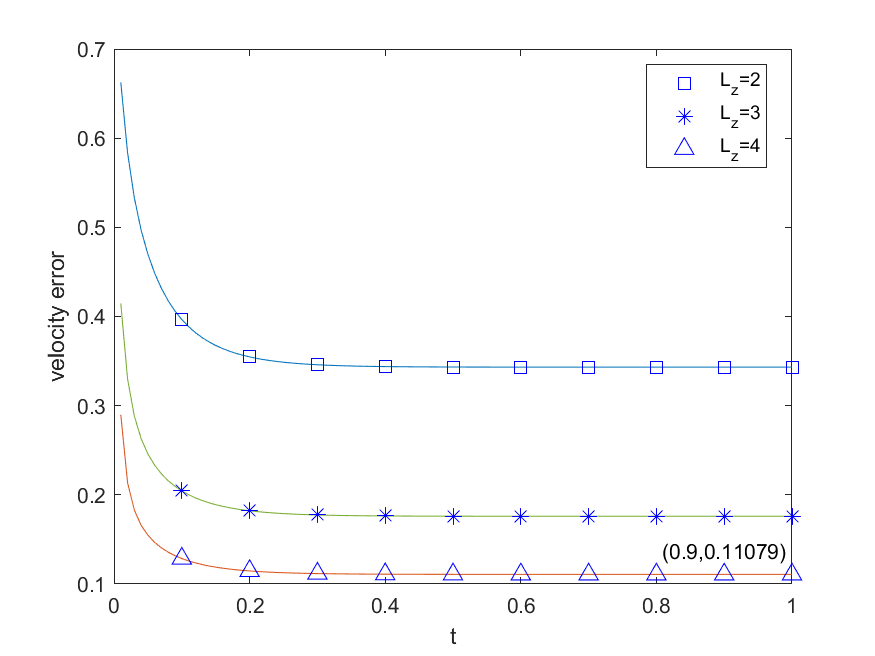}}
\subfigure[pressure error with $L_z$]
{\includegraphics[width=0.48\textwidth]{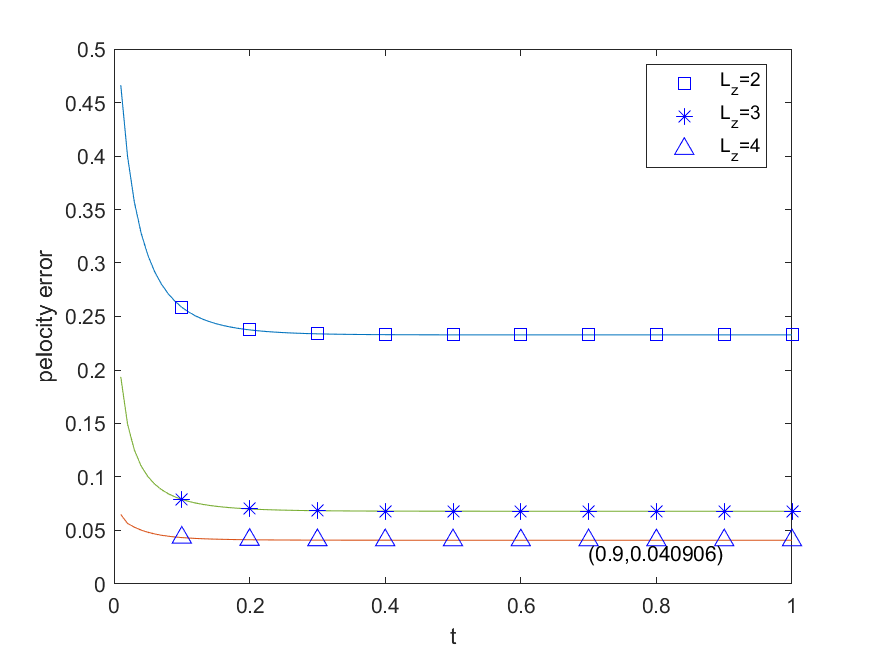}}
\subfigure[velocity error with $J$]
{\includegraphics[width=0.48\textwidth]{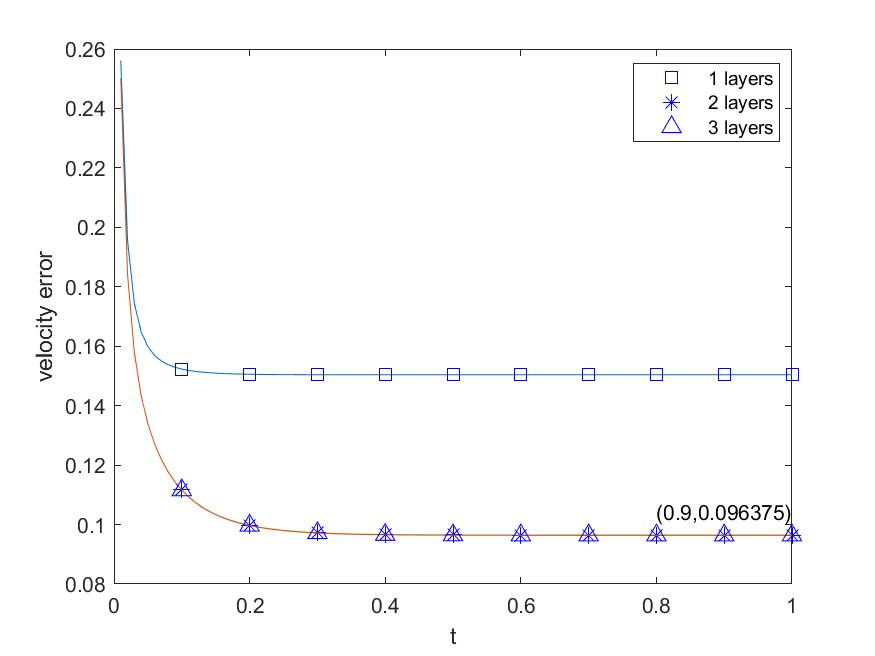}}
\subfigure[pressure error with $J$]
{\includegraphics[width=0.48\textwidth]{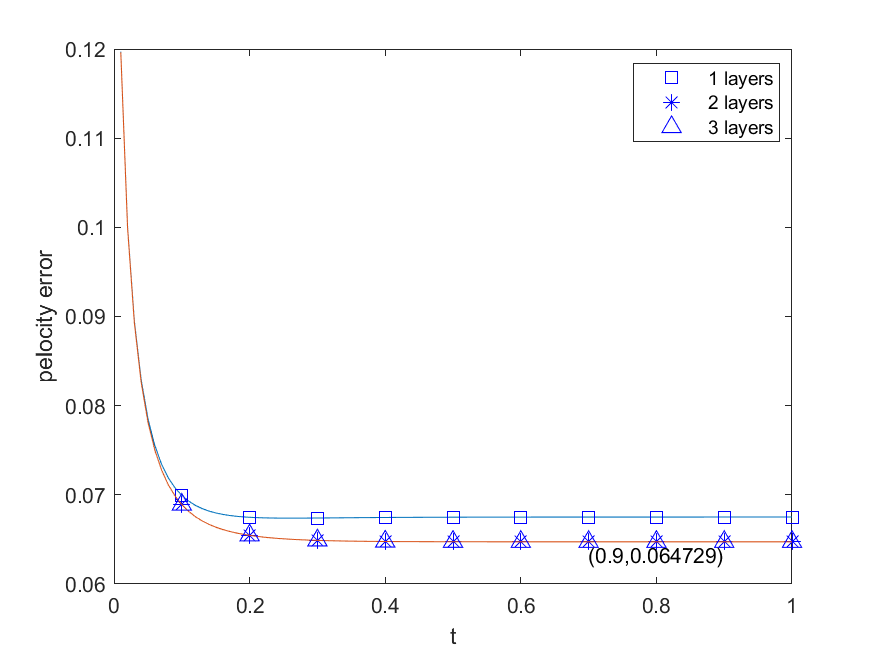}}
\subfigure[velocity error with $H$]
{\includegraphics[width=0.48\textwidth]{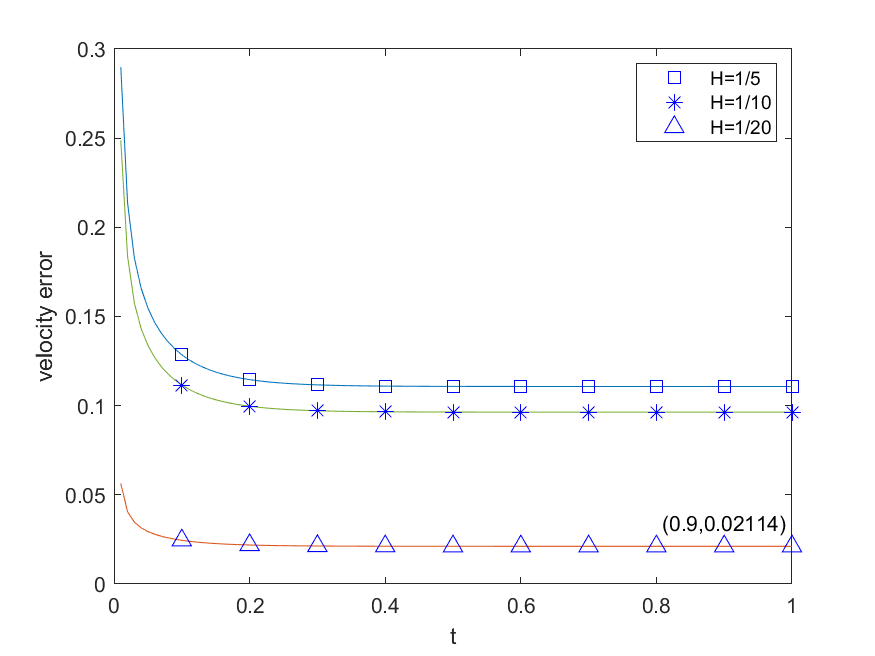}}
\subfigure[pressure error with $H$]
{\includegraphics[width=0.48\textwidth]{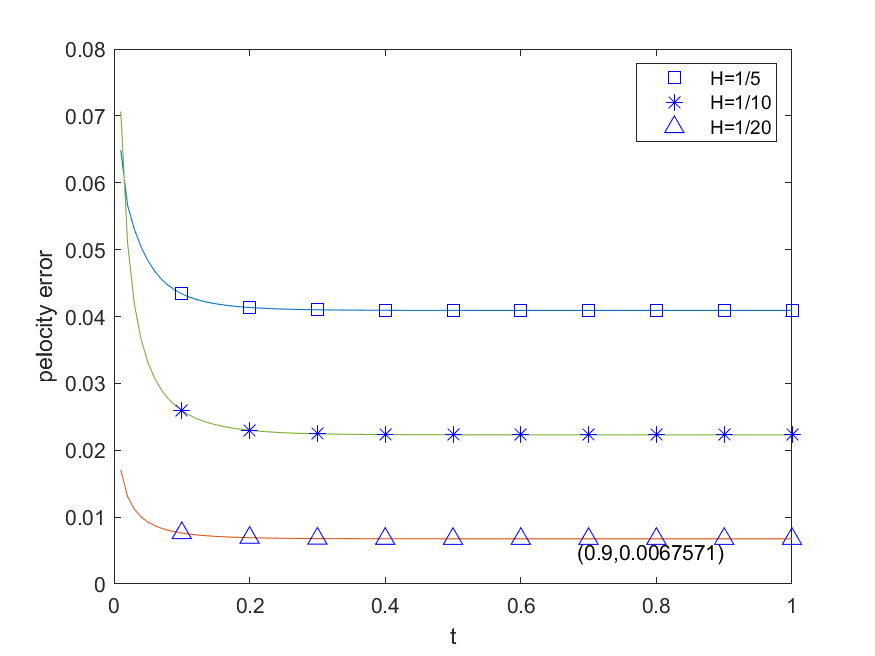}}
\caption{Example 2. Velocity errors and pressure errors with $\mathbf{\kappa_2}$ under different $L_z$, $J$ and $H$. First row : $J=2$, $H=1/5$. Second row: $L_z=2$, $H=1/10$. Third row: $L_z=4$, $J=4$.}\label{error_k2}
\end{figure}

In this part, we apply $\kappa_2$, which has the same meshsize $100\times 100$. Besides, we set terminal time $T=1$ and time step size $\tau=10^{-2}$ as in first case. In Figure \ref{vp_k2_t3}, we compare reference solutions and approximations obtained with three different settings of mixed CEM-GMsFEM. In this example, we compare three cases, $L_z=1,2,3$ and set $J=3$ with $H=1/10$. Similar to example 1, there are some apparent differences between the initial approximation $L_z=1$ and the reference. As for the last two cases ($L_z=3$ and $L_z=4$), they offer relatively better alternatives since there are fewer noticeable distinctions from the reference solution, which can be confirmed in Figure \ref{error_k2}. However, as one can observe from the figure, most of the complex fine information is remained in the last row.

In Figure \ref{error_k2}, we can see similar convergence results as in example 1. More specifically, in the first row, we compare three different $L_z$ under same $J=1$ and $H=1/5$, i.e. $L_z=2,3,4$. As one can observe from the graph, the error reduction is remarkable from over $30\%$ to about $10\%$ for velocity solutions. As for pressure solutions, there is a significant decay in error, which decreases from about $25\%$ to nearly $4\%$. In terms of the second row, we observe quick decay in error as one increase the oversampling layer from initial case. When $L_z=2$ and $H=1/10$, the velocity and pressure errors are decreased to about $9\%$ and $6\%$ respectively. For the last row, one can observe quicker convergence in the velocity solution, since results from the last two cases are almost identical. For the pressure solution, one can consistently see error reduction as we have finer coarse mesh. Consequently, one can get relatively good approximation with $L_z=4$, $J=4$ and $H=1/20$ with this permeability field.


\subsection{Example 3}
\begin{figure}[!htbp]
    \includegraphics[width=6in]{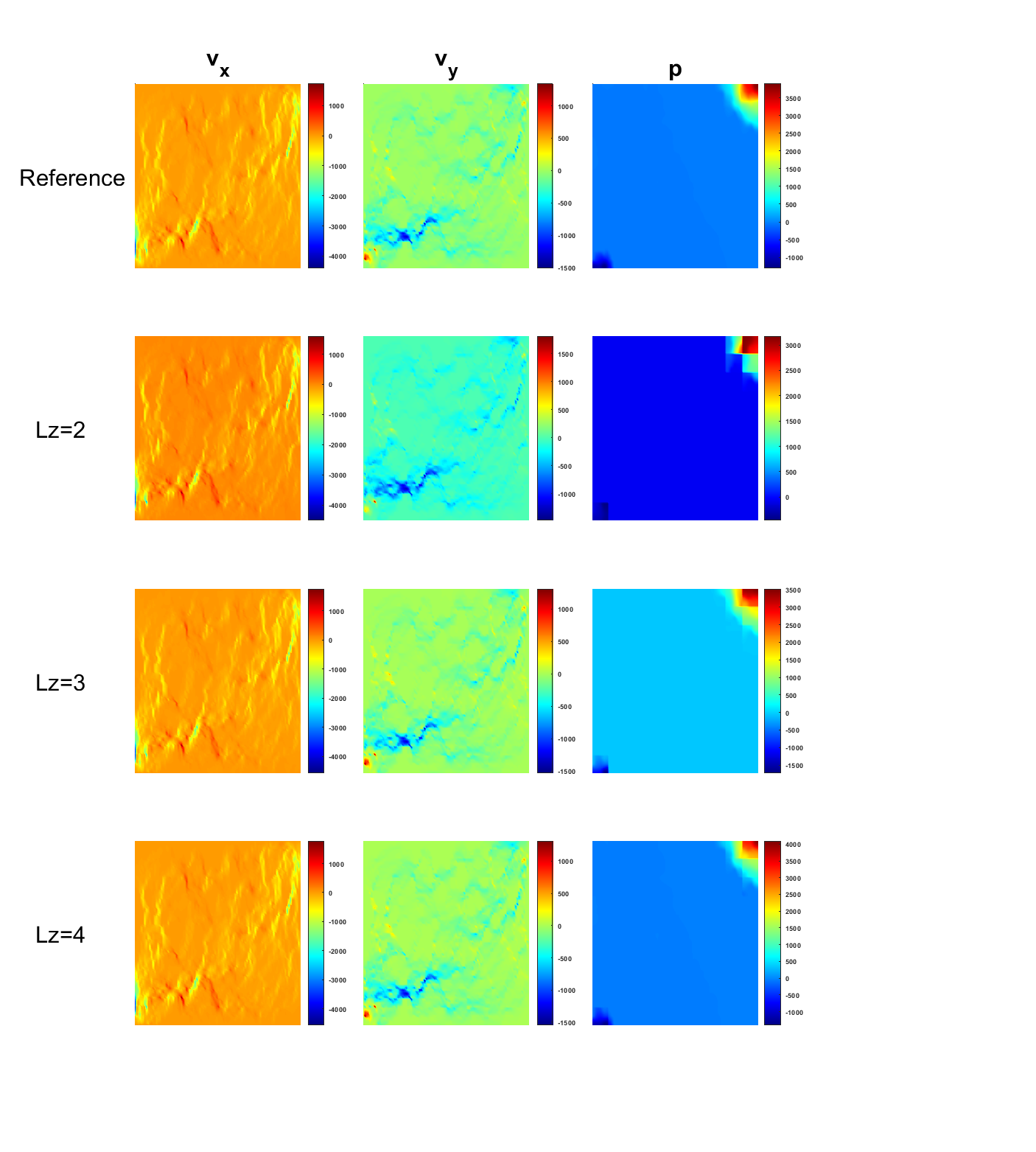}
    \caption{Example 3. Profile of velocity and pressure with $\mathbf{\kappa_3}$ at $T$ obtained by four methods. $J=1$ and $H=1/10$. Components of velocity solutions in x and y direction are shown in first two columns while pressure solution is in last column. Results from MFEM (reference), mixed CEM-GMsFEM  with $L_z=2$, $L_z=3$ and $L_z=4$ are shown in the first to fourth rows. }
    \label{vp_k3_t3}
\end{figure}

In this example, we apply $\kappa_3$ with fine meshsize $100\times 100$. Besides, we set terminal time $T=1$ and time step size $\tau=10^{-2}$. As we mentioned before, this field is most heterogeneous, which can be verified in Figure \ref{medium}. Different from previous two examples, as we can see from Figure \ref{vp_k3_t3}, there are apparent and steady improvements in the approximations as one increase the dimension of the multiscale space. More specifically, in initial case where $L_z=2$, some noticeable differences can be detected from three plots $v_x$, $v_y$ and $p$. For example, some fine-scare information in $v_y$ near the lower-left corner are neglected as well as the mismatched range of $v_y$. The approximation is improved in the $L_z=3$ case but some local information is distinct from reference solution like the upper-right corner in pressure solution. The last case $L_z=4$ offers a best alternative by keeping most fine-scale information.

\begin{figure}[!htbp]
\subfigure[velocity error with $L_z$]
{\includegraphics[width=0.48\textwidth]{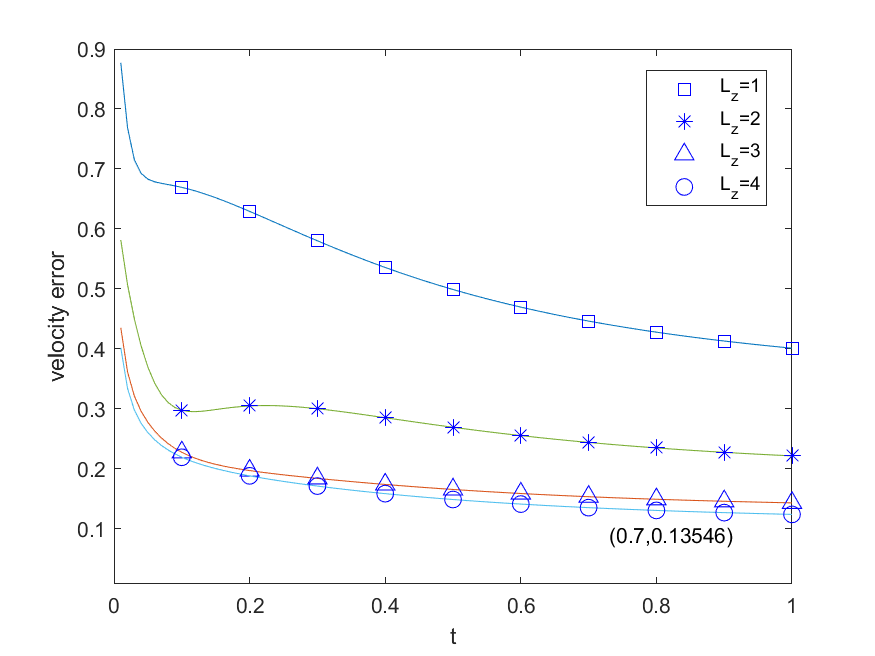}}
\subfigure[pressure error with $L_z$]
{\includegraphics[width=0.48\textwidth]{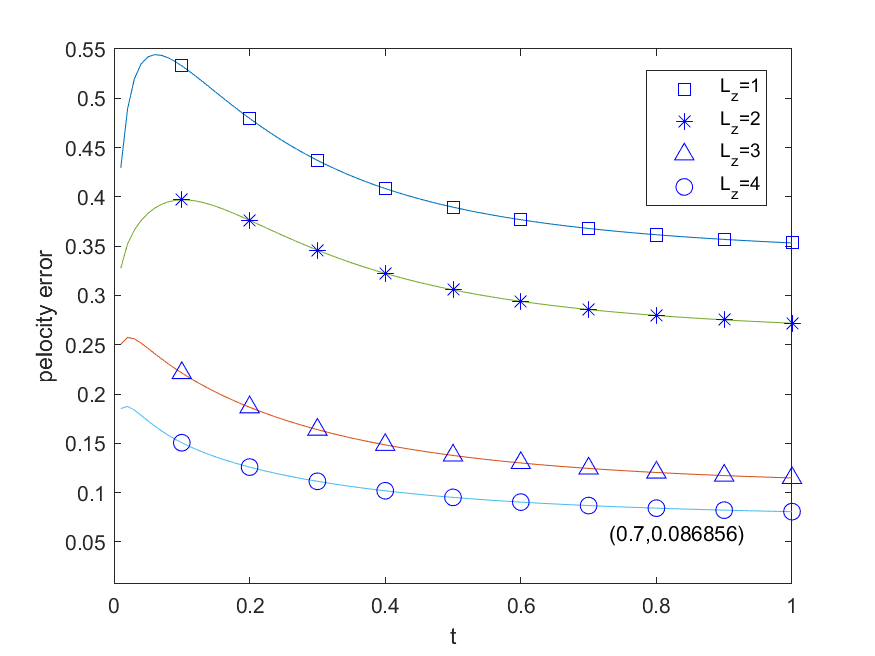}}
\subfigure[velocity error with $J$]
{\includegraphics[width=0.48\textwidth]{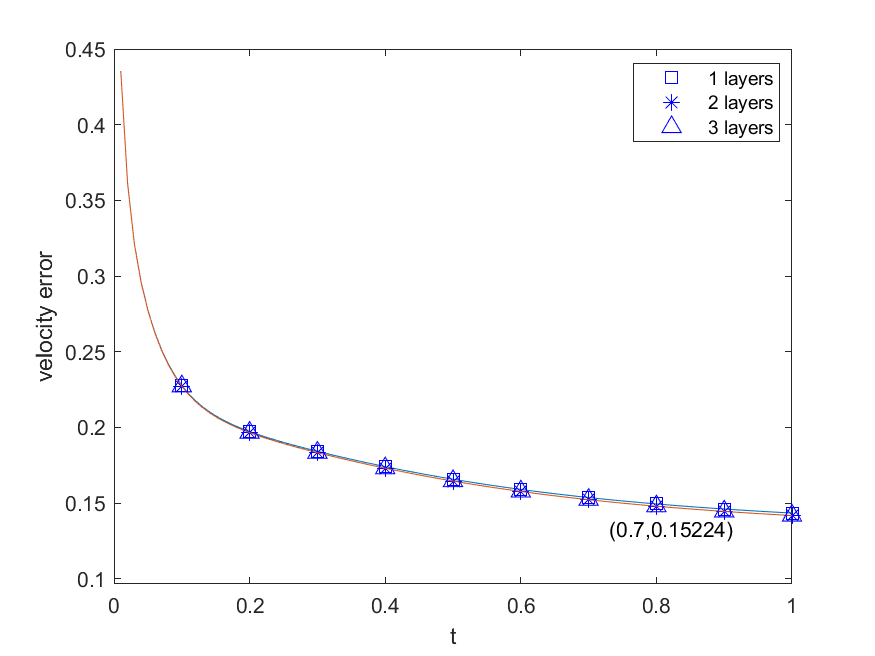}}
\subfigure[pressure error with $J$]
{\includegraphics[width=0.48\textwidth]{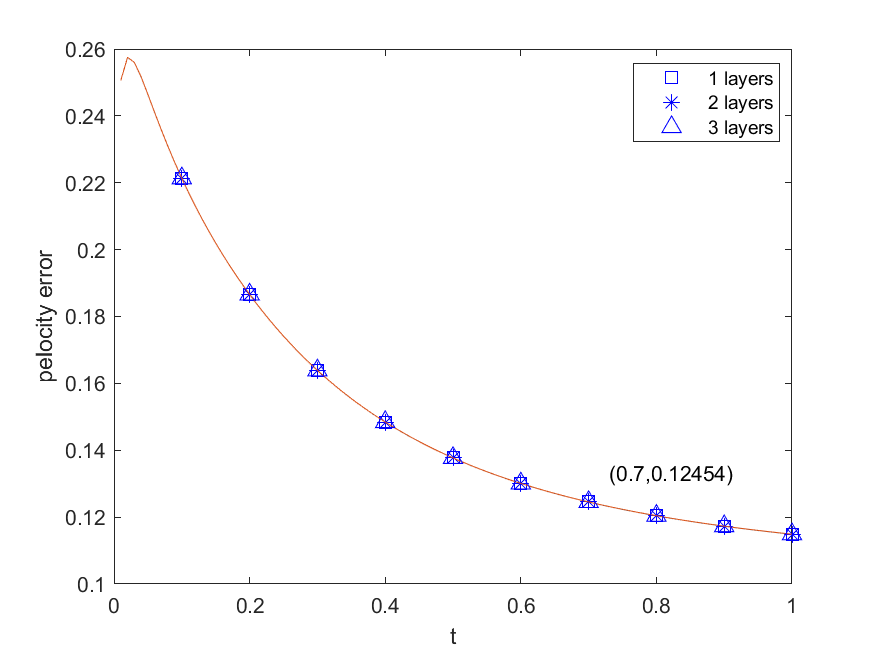}}
\subfigure[velocity error with $H$]
{\includegraphics[width=0.48\textwidth]{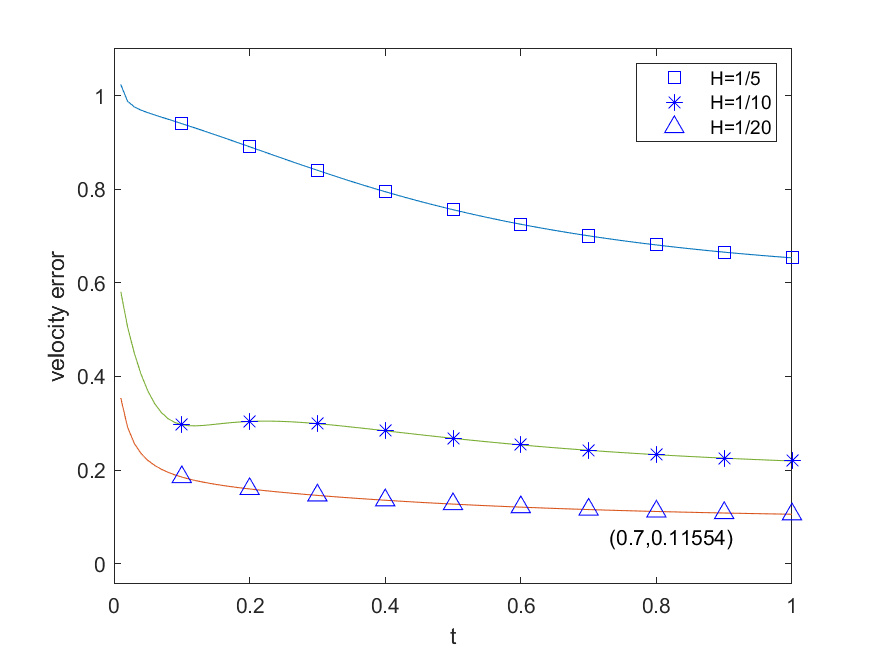}}
\subfigure[pressure error with $H$]
{\includegraphics[width=0.48\textwidth]{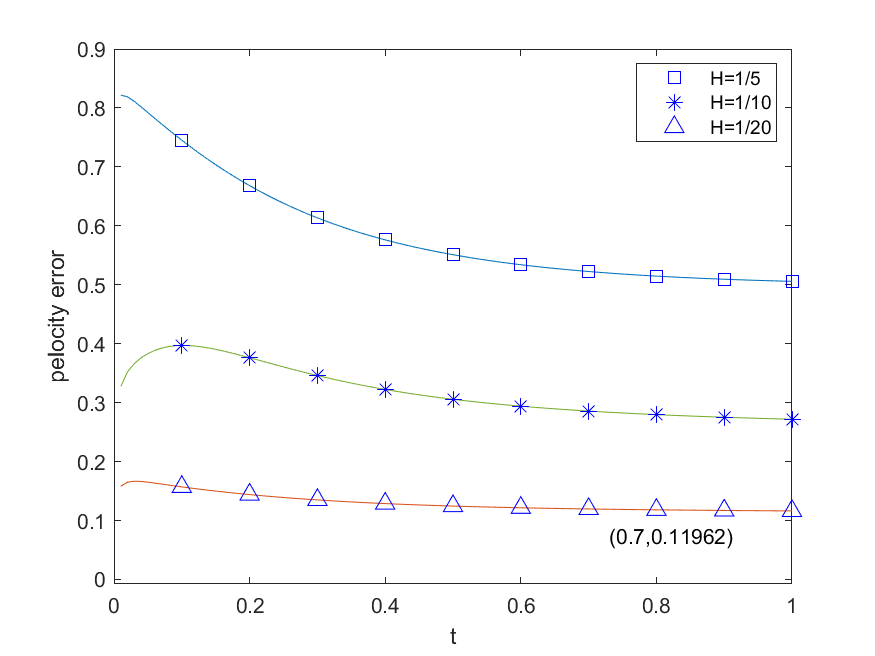}}
\caption{Example 3. Velocity errors and pressure errors with $\mathbf{\kappa_3}$ under different $L_z$, $J$ and $H$. First row : $J=1$, $H=1/10$. Second row: $L_z=3$, $H=1/10$. Third row: $L_z=2$, $J=2$.}\label{error_k3}
\end{figure}

In first row of Figure \ref{error_k3}, we compare four cases, from $L_z=1$ to $L_z=4$ and we set $J=1$ with $H=1/10$. One can see convergence for both velocity error and pressure error as increasing number of basis functions. For the former one, it is reduced from near $50\%$ to about $13\%$ while for the pressure error, the error decay is more remarkable, from about $40\%$ to $8\%.$ However, the improvement from $L_z=3$ to $L_z=4$ is less significant, which is similar to the former examples. For the second row, we set $L_z=3$ and $H=1/10$ and change $J$. As one may refer from the figure, one layer is sufficient in this case, and no further evident increase in accuracy can be observed by adding more layers, which is similar for both velocity and pressure solutions. As for the third row, we set $L_z=2$ and $J=2$, and we can see good convergence as one decrease $H$. Distinct from the second example, each accuracy improvement resulted from halving the coarse elements is remarkable. For the velocity solution, the error is diminished by half when $H$ is decreased to $1/20$ from $1/10$, which again verify the convergence analysis. And the error reduction is even more significant in pressure solution. Above all, one may easily see the efficiency of mixed CEM-GMsFEM in error reduction in this case.
\subsection{Example 4}
\begin{figure}[!htbp]
    \includegraphics[width=6in]{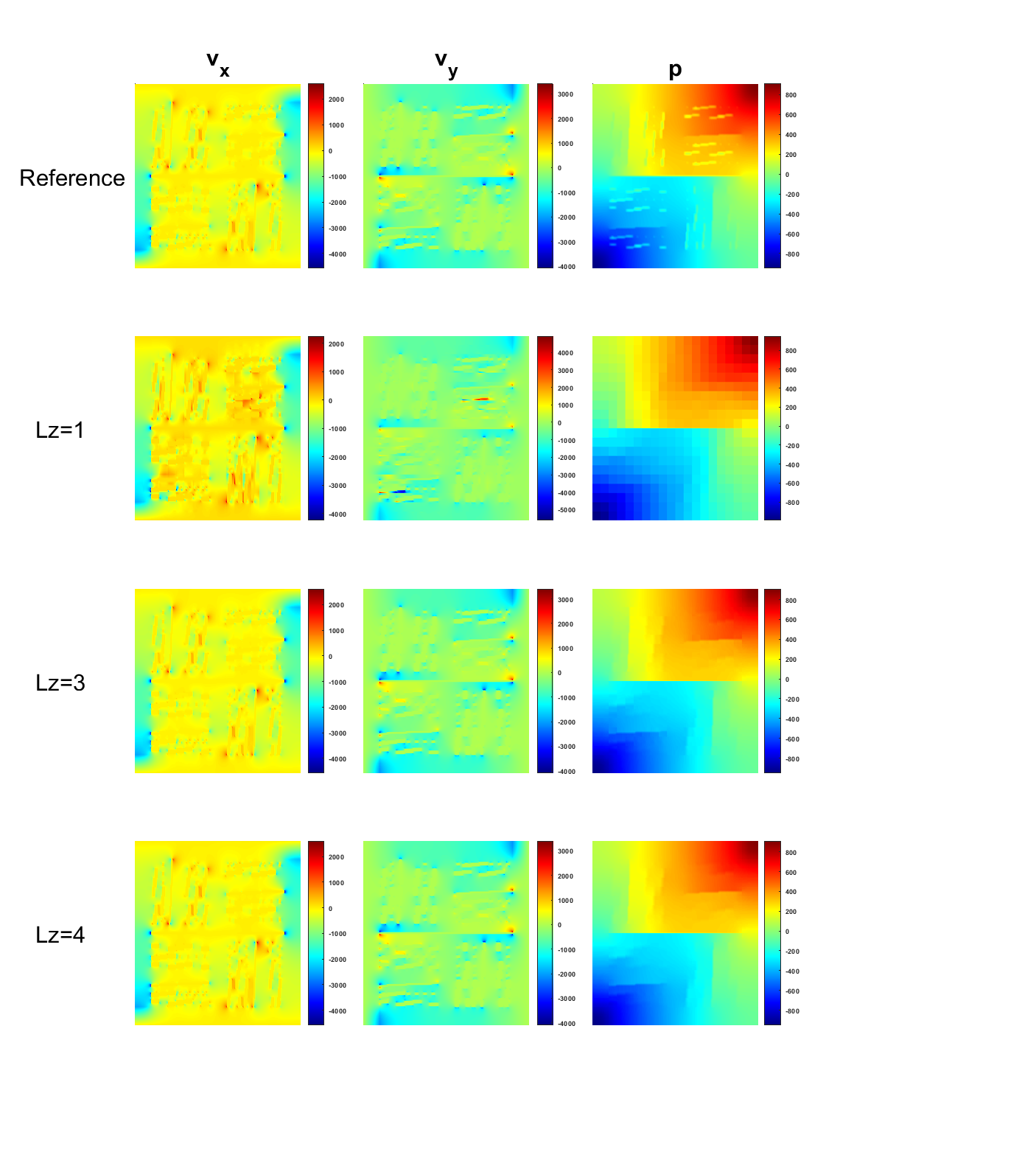}\caption{Example 4. Profile of velocity and pressure with $\mathbf{\kappa_4}$ at $T$ obtained by four methods. $J=2$ and $H=1/20$. Components of velocity solutions in x and y direction are shown in first two columns while pressure solution is in last column. Results from MFEM (reference), mixed CEM-GMsFEM  with $L_z=1$, $L_z=3$ and $L_z=4$ are shown in one to fourth rows.  }
    \label{vp_k4_t3}
\end{figure}

In this example, we consider $\kappa_4$ with fine size $200*200$. We also set the terminal time $T=1$ and time step size $\tau=10^{-2}$. Since $\kappa_4=1/\kappa_1$, these two share some similarities. For example the locations of the transitions are the same. However the locations with largest value corresponding to that with smallest value. In terms of other places, the permeability remains the same. In Figure \ref{vp_k4_t3}, we use $J=2$ and $H=1/20$ and plot the velocity and pressure at terminal time with different $L_z$. In the initial case where $L_z=1$, there are apparent differences from the reference cases. In particular, we can see some discontinuity on local neighborhood boundaries for pressure while the smoothness improves in the latter two cases. However all three approximation leave out some fine-scale information.

\begin{figure}[!htbp]
\subfigure[velocity error with $L_z$]
{\includegraphics[width=0.48\textwidth]{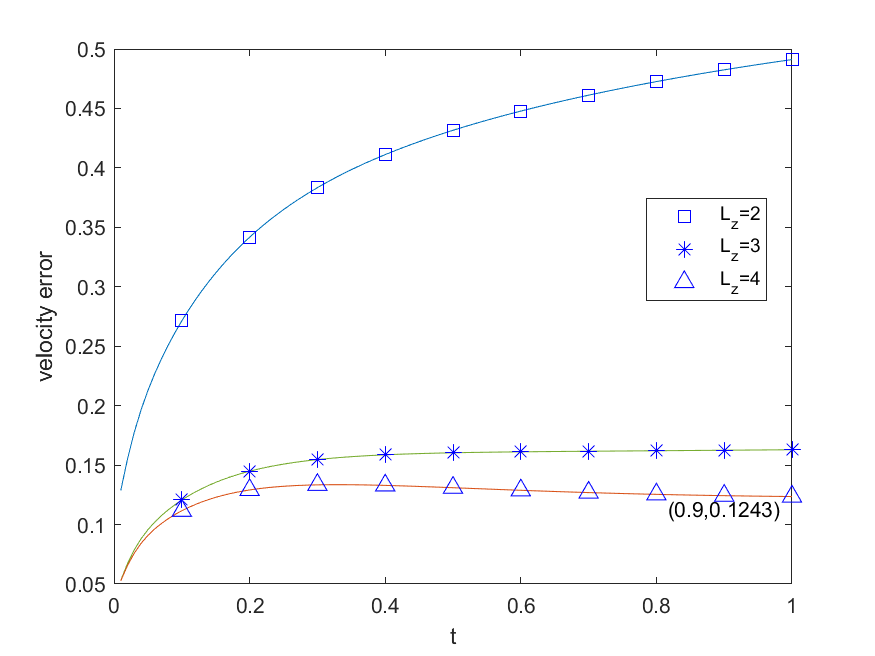}}
\subfigure[pressure error with $L_z$]
{\includegraphics[width=0.48\textwidth]{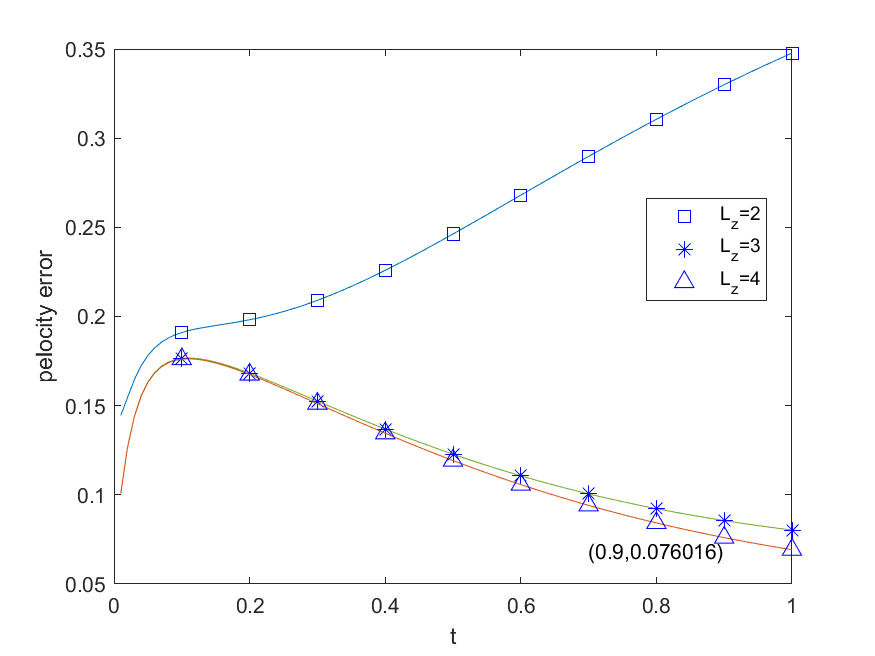}}
\subfigure[velocity error with $J$]
{\includegraphics[width=0.48\textwidth]{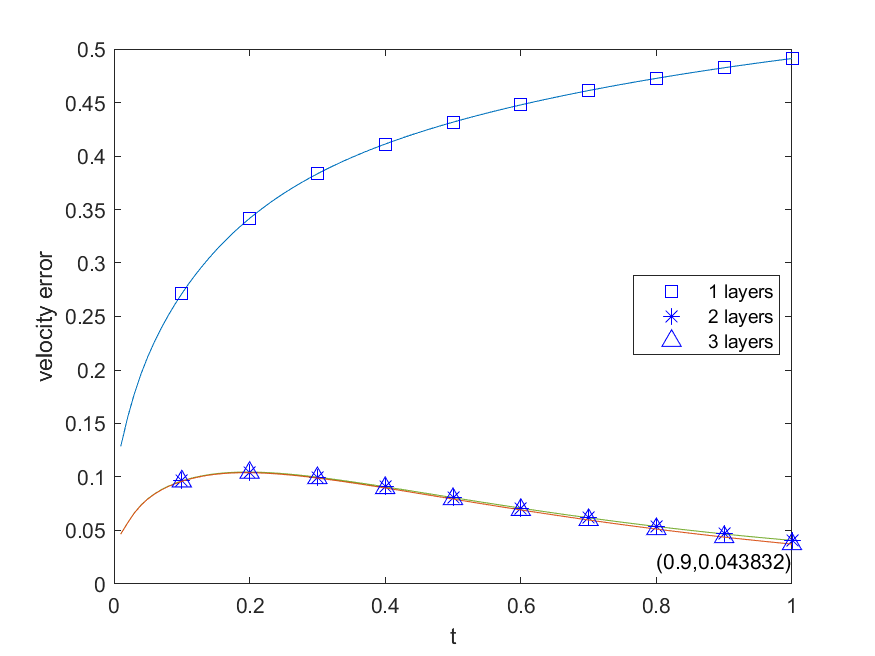}}
\subfigure[pressure error with $J$]
{\includegraphics[width=0.48\textwidth]{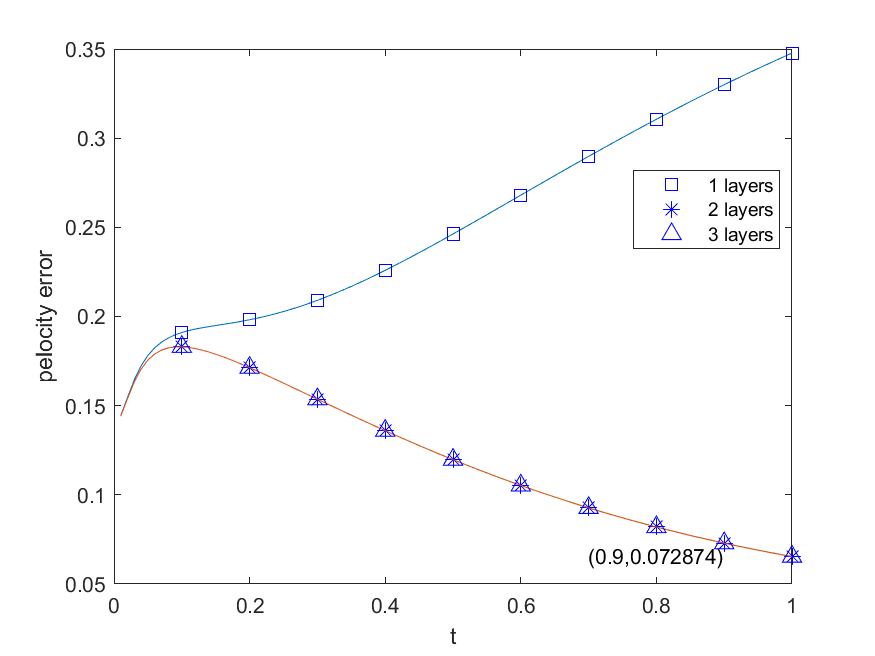}}
\subfigure[velocity error with $H$]
{\includegraphics[width=0.48\textwidth]{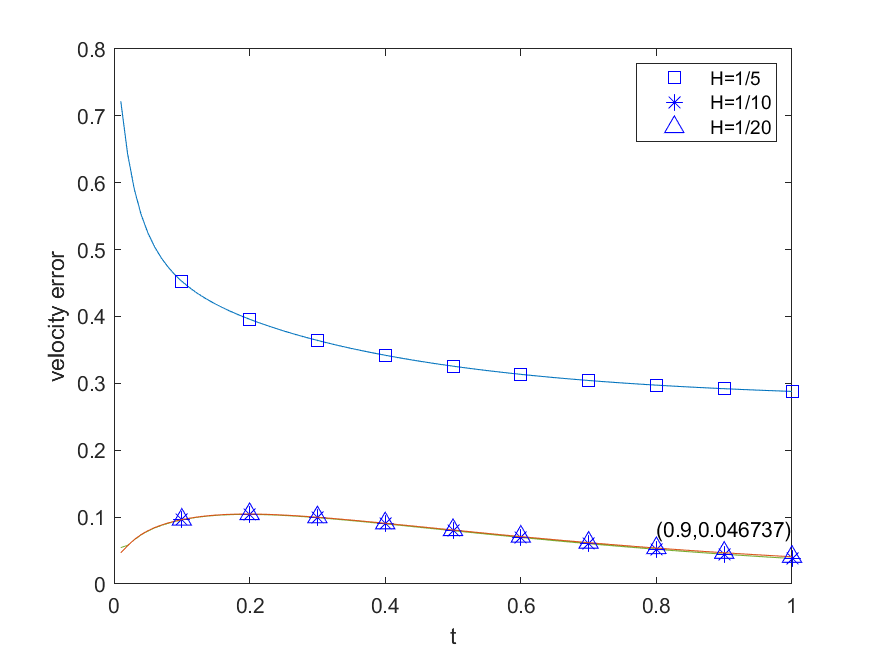}}
\subfigure[pressure error with $H$]
{\includegraphics[width=0.48\textwidth]{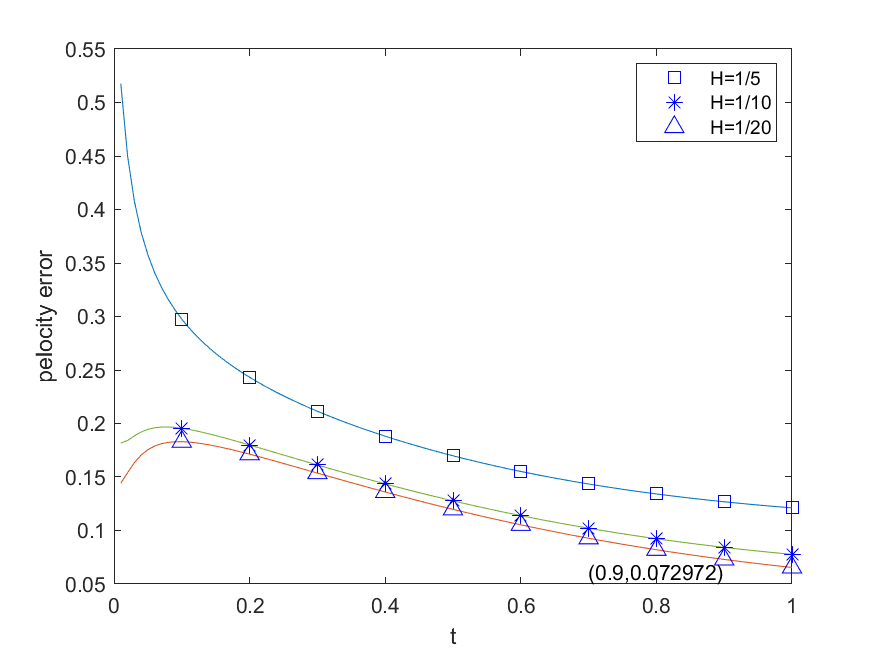}}
\caption{Example 4. Velocity errors and pressure errors with $\mathbf{\kappa_2}$ under different $L_z$, $J$ and $H$. First row : $J=1$, $H=1/20$. Second row: $L_z=2$, $H=1/20$. Third row: $L_z=2$, $J=2$.}\label{error_k4}
\end{figure}

In Figure \ref{error_k4}, we still consider the velocity and pressure error under different settings. In the first row, we consider use different $L_z$ with $J=1$ and $H=1/20$. We can see the in the first case with only one coarse oversampling layer, the error can not decrease as time, which is true in other two cases. In the second row, we can observe similar effects and the improvements are resulted by increasing the number of oversampling layers. From here, we may see that one can improve the accuracy efficiently either by using bigger local neighborhood or increasing the basis used in the local neighborhood. From the third row, one may see the accuracy is significantly improved by using finer coarse neighborhoods. As when we double the number of coarse neighborhoods starting at $H=1/5$, the velocity error is largely decreased from over $30\%$ to about $4\%$. Last but not least, one may see the importance of using sufficiently large oversampled coarse neighborhood to obtain relatively good approximations, which correspond to remarks by the end of Section \ref{analysis}.

\section{Conclusion}
In this paper, we solve heterogeneous parabolic equation with Constraint energy minimization generalized multiscale finite element method in mixed formulation. We review the construction of the multiscale space, where the key part is solving some well-designed local spectral problems. Furthermore, we prove the convergence of the proposed method. In the numerical results, we test the proposed method with four distinctive permeability fields. We mainly compare the influence of three parameters in the multiscale methods, which are $L_z$, $J$, $H$. Our results  can well verify the convergence analysis.

\section*{Acknowledgments}

The research of Eric Chung is partially supported by the Hong Kong RGC General Research Fund (Project numbers 14304719 and 14302018) and the CUHK Faculty of Science Direct Grant 2019-20.

\bibliographystyle{plain}
\bibliography{reference}
\end{document}